\documentclass[a4paper,10pt,twoside]{scrartcl}
\usepackage{german,amssymb,amsmath,latexsym,epsfig,graphics,graphicx,mathrsfs,color,amsthm,cite}
\usepackage[latin1]{inputenc}
\usepackage{makeidx}

\newtheorem{theorem}{Theorem}[section]
\newtheorem{remark}[theorem]{Remark}
\newtheorem{corollary}[theorem]{Corollary}
\newtheorem{definition}[theorem]{Definition}
\newtheorem{lemma}[theorem]{Lemma}
\newtheorem{example}[theorem]{Example}

\numberwithin{equation}{section}

\newcommand{\R}{{\mathbb R}}
\newcommand{\N}{{\mathbb N}}

\DeclareMathOperator*{\esssup}{ess\,sup}

\begin{document}
\title{\Large{Pontryagin's Maximum Principle for Infinite Horizon Optimal Control Problems with Bounded Processes and
       with State Constraints}}

\author{\large{Nico Tauchnitz} \\ {\small Nico.Tauchnitz@hszg.de}}
\date{\large{\today}}
\maketitle

\begin{abstract}
\textbf{Abstract.} This paper concerns a class of infinite horizon optimal control problems with state constraints.
By extending the needle variation method to the infinite horizon case
we obtain a complete set of necessary optimality conditions for a strong local minimizer in
form of the Pontryagin maximum principle and
the validity of several transversality conditions at infinity.
To this purpose, we develop an approach in the space of continuous functions converging at infinity.
The considerations including new results on optimality and transversality conditions,
on Banach spaces of continuous functions,
on the needle variation method and on problems with state constraints on the unbounded time horizon. \\[2mm]
\textbf{Keywords.} Pontryagin Maximum Principle $\cdot$ Optimal Control $\cdot$ Infinite Horizon $\cdot$
                   Optimality Conditions $\cdot$ Transversality Conditions $\cdot$ State Constraints
\end{abstract}


\section{Introduction}
Infinite horizon optimal control is a particular class in control theory.
The origin of this  class goes back to the work of Ramsey \cite{Ramsey} in 1928,
where the optimization of economic growth is considered as a variational problem on an unbounded planning horizon.
With the discovery of the Pontryagin maximum principle for standard optimal control problems by Pontryagin et. al. \cite{Pontrjagin} in 1958, 
the impetus for the derivation of necessary conditions in the infinite horizon case was initiated.
In the recent years,
new results were established by different approaches like the finite horizon approximation, the needle variation technique,
the approach by extremal principles, the discrete time framework or the dynamic programming
(cf. \cite{AseKry,AseVel,AseVel2,AseVel3,Pickenhain,LyPi,Blot,BachirBlot,Frankowska,Frankowska2} and the references therein). \\
In this paper,
we present the Pontryagin maximum principle for infinite horizon optimal control problems with state constraints.
The proof bases on the procedure in Ioffe \& Tichomirov \cite{Ioffe} for regularly locally convex problems.
But this method must be extended to the infinite horizon case.
We achieve necessary conditions in form of the adjoint equation, the maximum condition and standard type transversality conditions.
Our framework is closely related to the paper of Brodskii \cite{Brodskii}.
But in contrast to \cite{Brodskii},
we consider strong local optimality and we use the state space of continuous functions converging at infinity.
The last fact delivers the complete information on the subdifferential of state constraints on the infinite horizon.
Consequently, we obtain a more aesthetic representation of the adjoint as in \cite{Brodskii}.
Moreover,
we investigate objectives with general summable functions like a deterministic, unbounded Weibull distribution. \\
The paper is organized as follows. 
We start with the statement of the control problem, the main assumptions and the maximum principle (Section \ref{SectionProblem}).
In Section \ref{SectionProof}, we present the proof of the maximum principle.
The proof is divided into different parts: the introduction of the spaces of continuous functions converging at infinity,
the statement and properties of the extremal problem,
results on needle variations and linear differential equations in the infinite horizon case,
the convex analysis of state constraints on the infinite horizon
and finally the completion of the proof by standard arguments.
Section \ref{SectionTransversality} is devoted to transversality conditions and to the normal form of the maximum principle
in the absence of state constraints.
Then we show Arrow type sufficiency conditions in Section \ref{SectionArrow}.
We conclude the paper with Section \ref{SectionRelation},
where we show the relation between finite horizon and infinite horizon optimal control problems.
\section{Problem Formulation and the Maximum Principle} \label{SectionProblem}
We consider the following infinite horizon optimal control problem:
\begin{eqnarray}
&& \label{Aufgabe1} J\big(x(\cdot),u(\cdot)\big) = \int_0^\infty \omega(t) f\big(t,x(t),u(t)\big) \, dt \to \inf, \\
&& \label{Aufgabe2} \dot{x}(t) = \varphi\big(t,x(t),u(t)\big), \\
&& \label{Aufgabe3} h_0\big(x(0)\big)=0, \qquad \lim_{t \to \infty} h_1\big(t,x(t)\big)=0,\\
&& \label{Aufgabe4} u(t) \in U \subseteq \R^m, \quad U \not= \emptyset, \\
&& \label{Aufgabe5} g_j\big(t,x(t)\big) \leq 0, \quad t \in \R_+, \quad j=1,...,l.
\end{eqnarray}
Throughout this paper let $\R_+=[0,\infty)$ and $\overline{\R}_+ = [0,\infty]$.
In the problem (\ref{Aufgabe1})--(\ref{Aufgabe5}) let the density function $\omega(\cdot)$ be an element of the space $L_1(\R_+,\R)$ with
$\|\omega(\cdot)\|_{L_1} >0$
and let $f:\R \times \R^n \times \R^m \to \R$, $\varphi:\R \times \R^n \times \R^m \to \R^n$, $g_j:\R \times \R^n \to \R$, $j=1,...,l$,
$h_0:\R^n \to \R^{s_0}$, $h_1:\R \times \R^n \to \R^{s_1}$. \\[2mm]
A trajectory $x(\cdot)$ of system (\ref{Aufgabe2}) corresponding to a control $u(\cdot) \in L_\infty(\R_+,U)$ is a solution
(on any bounded interval in the sense of Carath\'eodory \cite{Filippov})
of system (\ref{Aufgabe2}) that is defined on $[0,\infty)$.

\begin{definition}
The set of all admissible processes $\mathscr{A}_{adm}$ consists of all
$\big(x(\cdot),u(\cdot)\big) \in W^1_\infty(\R_+,\R^n) \times L_\infty(\R_+,U)$
satisfying (\ref{Aufgabe2})--(\ref{Aufgabe5}) and making the Lebesgue integral in (\ref{Aufgabe1}) finite.
\end{definition}

Let $x(\cdot)\in W^1_\infty(\R_+,\R^n)$ and
$V_\gamma= \{ (t,x) \in \overline{\R}_+ \times \R^n\,|\, \|x-x(t)\|\leq \gamma\}$ with $\gamma>0$.
Then $\mathscr{X}_{Lip}$ consists of all $x(\cdot)\in W^1_\infty(\R_+,\R^n)$ with the following properties:
For any compact subset $U_1$ of $\R^m$ there exists a number $\gamma >0$ such that
the mappings $f(t,x,u)$, $\varphi(t,x,u)$, $h_0(x)$, $h_1(t,x)$ and $g_j(t,x)$, $j=1,...,l$, are uniformly continuous w.r.t. all variables and
uniformly continuously differentiable w.r.t. $x$ on $V_\gamma \times U_1$.

\begin{definition}
The set $\mathscr{A}_{Lip}$ denotes the set of all pairs
$\big(x(\cdot),u(\cdot)\big) \in W^1_\infty(\R_+,\R^n) \times L_\infty(\R_+,U)$ with
$x(\cdot) \in \mathscr{X}_{Lip}$.
\end{definition}
According to the special character of infinite horizon optimal control problems with bounded processes
we state the following restrictive properties:
Let $\big(x(\cdot),u(\cdot)\big) \in W^1_\infty(\R_+,\R^n) \times L_\infty(\R_+,U)$ and
\begin{equation} \label{condition1}
\int_0^\infty \big\|\varphi\big(t,x(t),u(t)\big)\big\| \, dt < \infty, \qquad \int_0^\infty \big\|\varphi_x\big(t,x(t),u(t)\big)\big\| \, dt < \infty.
\end{equation}
Furthermore, we assume that for any $\delta>0$ there exists a $T>0$ such that the relation
\begin{equation} \label{condition2} 
\int_T^\infty \big\| \varphi\big(t,\xi(t),u(t)\big)-\varphi\big(t,\xi'(t),u(t)\big) - \varphi_x\big(t,x(t),u(t)\big)\big(\xi(t)-\xi'(t)\big) \big\| \, dt
\leq \delta \|\xi(\cdot)-\xi'(\cdot)\|_\infty
\end{equation}
holds for all $\xi(\cdot), \xi'(\cdot) \in W^1_\infty(\R_+,\R^n)$ with $\|\xi(\cdot)-x(\cdot)\|_\infty \leq \gamma$, $\|\xi'(\cdot)-x(\cdot)\|_\infty \leq \gamma$.

\begin{definition}
The set $\mathscr{A}_{lim}$ denotes the set of all pairs
$\big(x(\cdot),u(\cdot)\big) \in W^1_\infty(\R_+,\R^n) \times L_\infty(\R_+,U)$
satisfying (\ref{condition1}) and (\ref{condition2}).
\end{definition}

\begin{remark}
According to the first condition in (\ref{condition1})
for any $\big(x(\cdot),u(\cdot)\big) \in \mathscr{A}_{adm} \cap \mathscr{A}_{lim}$ the state trajectory $x(\cdot)$
possesses a limit at infinity.
\end{remark}

\begin{definition}
An admissible process $\big(x_*(\cdot),u_*(\cdot)\big)$ is a strong local minimizer
iff there exists a number $\varepsilon > 0$ such that the inequality
$J\big(x(\cdot),u(\cdot)\big) \geq J\big(x_*(\cdot),u_*(\cdot)\big)$
holds for any $\big(x(\cdot),u(\cdot)\big) \in \mathscr{A}_{adm}$
for which $\| x(\cdot)-x_*(\cdot) \|_\infty \leq \varepsilon$.
\end{definition}

\begin{definition} \label{DefinitionBorel}
\begin{enumerate}
\item[(a)] The Borel algebra on $\overline{\R}_+$ consists of all sets $B=A \cup E$ with open subsets $A$ of $\R_+$ and
           $E \subseteq \{\infty\}$.
\item[(b)] By $\mathscr{M}(\overline{\R}_+)$ we denote the set of the Borel measures on $\overline{\R}_+$,
           which have a unique representation $\mu=\mu_0+\mu_\infty$,
           where $\mu_0 $ is a regular signed Borel measure on $\R_+$ (cf. \cite{Rudin})
           and $\mu_\infty=\mu(\{\infty\})$ is a finite signed measure concentrated at $t=\infty$.
\end{enumerate}
\end{definition}

We introduce the Pontryagin function $H: \R \times \R^n \times \R^m \times \R^n \times \R \to \R$,
$$H(t,x,u,p,\lambda_0) = -\lambda_0 \omega(t) f(t,x,u) + \langle p, \varphi(t,x,u) \rangle.$$

\begin{theorem} \label{SatzPMPZB} 
Let $\big(x_*(\cdot),u_*(\cdot)\big) \in \mathscr{A}_{adm} \cap \mathscr{A}_{Lip} \cap \mathscr{A}_{lim}$. 
If $\big(x_*(\cdot),u_*(\cdot)\big)$ is a strong local minimizer in the problem (\ref{Aufgabe1})--(\ref{Aufgabe5}),
then there exist a number $\lambda_0 \geq 0$, vectors $l_0 \in \R^{s_0}$ and $l_1 \in \R^{s_1}$,
a vector-valued function $p(\cdot):\R_+ \to \R^n$ and non-negative Borel measures $\mu_j \in \mathscr{M}(\overline{\R}_+)$,
$j=1,...,l$, supported on the sets
$T_j=\big\{ t \in \overline{\R}_+ \,|\, g_j\big(t,x_*(t)\big)=0\big\}$, respectively,
not all zero and such that
\begin{enumerate}
\item[(a)] the vector-valued function $p(\cdot)$ is a solution of the integral equation
           \begin{equation}\label{PMPZB1}
           p(t) = - \lim_{t \to \infty}h_{1x}^T\big(t,x_*(t)\big) l_1 + \int_t^\infty H_x\big(s,x_*(s),u_*(s),p(s),\lambda_0\big) \, ds
                  - \sum_{j=1}^l \int_t^\infty g_{jx}\big(s,x_*(s)\big) \, d\mu_j(s)
           \end{equation}
           and
           \begin{equation}\label{PMPZB2}
           p(0) = {h_0'}^T\big(x_*(0)\big)l_0;
           \end{equation}
\item[(b)] for almost all $t \in \R_+$ the maximum condition holds:    
           \begin{equation}\label{PMPZB3}
           H\big(t,x_*(t),u_*(t),p(t),\lambda_0\big) = \max_{u \in U} H\big(t,x_*(t),u,p(t),\lambda_0\big).
           \end{equation}
\end{enumerate}
\end{theorem}

In the presence of state constraints, the adjoint $p(\cdot)$ can be discontinuous.
However, this function is of bounded variation and continuous from the left.
In the case, the measures $\mu_j$ containing non-zero masses concentrated at infinity, it follows
\begin{equation} \label{Transversality}
\lim_{t \to \infty} p(t)= \lim_{t \to \infty} \bigg[-h_{1x}^T\big(t,x_*(t)\big) l_1 - \sum_{j=1}^l g_{jx}\big(t,x_*(t)\big) \mu_j(\{\infty\}) \bigg].
\end{equation}

The requirement $\big(x_*(\cdot),u_*(\cdot)\big) \in \mathscr{A}_{\lim}$ is a major restriction for the application of Theorem \ref{SatzPMPZB}.
We demonstrate this fact in the following example.
\begin{example} \label{BeipielRWUnendlich1}
We consider the problem
\begin{eqnarray*}
&& J\big(x(\cdot),z(\cdot),u(\cdot)\big)=\int_0^\infty e^{-\varrho t}\big(1-u(t)\big)x(t) \, dt \to \sup, \\
{\mathit s.t.} && \dot{x}(t)=u(t)x(t), \qquad x(0)=1, \qquad u \in [0,1], \qquad \varrho \in (0,1) \\
&& \dot{z}(t)=e^{-\varrho t}x(t), \qquad z(0)=0, \qquad z(t) \leq Z,\; t \in \R_+, \qquad Z> \frac{1}{\varrho}.
\end{eqnarray*}
The latter differential equation and the state constraint are the result of the budget constraint
$$\int_0^\infty e^{-\varrho t}x(t) \, dt \leq Z.$$
Since $\dot{z}(t) >0$ on $\R_+$ and $z(t) \leq Z$ it follows $e^{-\varrho t}x(t) \to 0$ for $t \to \infty$.
Then the inequaility
\begin{eqnarray*}
    J\big(x(\cdot),z(\cdot),u(\cdot)\big)
&=& \int_0^\infty e^{-\varrho t}\big(1-u(t)\big)x(t) \, dt = \int_0^\infty \dot{z}(t) \, dt - \int_0^\infty e^{-\varrho t} \dot{x}(t) \, dt \\
&=& \int_0^\infty \dot{z}(t) \, dt + 1 - \varrho \int_0^\infty e^{-\varrho t} x(t) \, dt \leq 1 + (1-\varrho) Z
\end{eqnarray*}
shows,
that any admissible process with $\lim\limits_{t \to \infty} z(t) = Z$ is globally optimal.
We discuss two policies:
\begin{enumerate}
\item[(A)] A bang-bang type control delivers the admissible process $\big(x_*(\cdot),z_*(\cdot),u_*(\cdot)\big)$ with
           \begin{eqnarray*}
           x_*(t) &=& \left\{ \begin{array}{ll} e^t,& t \in [0,\tau), \\ e^\tau, & t \in [\tau,\infty), \end{array} \right. \quad
                      u_*(t) = \left\{ \begin{array}{ll} 1,& t \in [0,\tau), \\ 0, & t \in [\tau,\infty), \end{array} \right. \\
           z_*(t) &=& \left\{ \begin{array}{ll} \frac{1}{1-\varrho}\big(e^{(1-\varrho)t} - 1 \big) ,& t \in [0,\tau), \\
                                                z(\tau) + \frac{1}{\varrho}\big(e^{(1-\varrho)\tau} - e^{\tau-\varrho t} \big), & t \in [\tau,\infty).
             \end{array} \right. 
           \end{eqnarray*}
           Herein, the switching time $\tau >0$ satisfies the condition
           $$e^{(1-\varrho)\tau}\bigg(\frac{1}{\varrho}+\frac{1}{1-\varrho}\bigg) = Z+\frac{1}{1-\varrho}.$$
           Since $u_*(\cdot) \in L_1(\R_+,[0,1])$ the restrictive requirement $\big(x_*(\cdot),z_*(\cdot),u_*(\cdot)\big) \in \mathscr{A}_{\lim}$ holds.
           Therefore, the assumptions of Theorem \ref{SatzPMPZB} are satisfied.
           Applying the necessary conditions, we obtain the multipliers
           $$\lambda_0=1, \qquad p(t)= e^{-\varrho t}, \qquad q(t)= \varrho-1,$$
           where $q(\cdot)$ satisfies the transversality condition (\ref{Transversality}) with $\mu(\{\infty\})=1-\varrho$.
\item[(B)] A constant investment rate leads to the admissible process
           $\big(\xi_*(\cdot),\zeta_*(\cdot),w_*(\cdot)\big)$ with
           $$\xi_*(t) = e^{\alpha t},\qquad w_*(t)=\alpha,\qquad 
             \zeta_*(t)= \frac{1}{\alpha - \varrho} (e^{(\alpha-\varrho)t}-1),\qquad \alpha=\varrho-\frac{1}{Z} \in (0, \varrho).$$
           Since $w_*(\cdot) \not\in L_1(\R_+,[0,1])$, the requirements of Theorem \ref{SatzPMPZB} are violated. 
\end{enumerate}
\end{example}

In the absence of state constraints Theorem \ref{SatzPMPZB} leads to the following version of the maximum principle:

\begin{theorem} \label{SatzPMP} 
Let $\big(x_*(\cdot),u_*(\cdot)\big) \in \mathscr{A}_{adm} \cap \mathscr{A}_{Lip} \cap \mathscr{A}_{lim}$. 
If $\big(x_*(\cdot),u_*(\cdot)\big)$ is a strong local minimizer in the problem (\ref{Aufgabe1})--(\ref{Aufgabe4}),
then there exist a number $\lambda_0 \geq 0$, vectors $l_0 \in \R^{s_0}$, $l_1 \in \R^{s_1}$ and a vector-valued function $p(\cdot):\R_+ \to \R^n$,
not all zero, such that
\begin{enumerate}
\item[(a)] the vector-valued function $p(\cdot)$ satisfies almost everywhere on $\R_+$ the adjoint equation
           \begin{equation}\label{PMP1}
           \dot{p}(t) = -\varphi_x^T\big(t,x_*(t),u_*(t)\big) p(t) + \lambda_0 \omega(t)f_x\big(t,x_*(t),u_*(t)\big)
           \end{equation}
           and the transversality conditions
           \begin{equation}\label{PMP2}
           p(0) = {h_0'}^T\big(x_*(0)\big)l_0, \qquad \lim_{t \to \infty} p(t)= - \lim_{t \to \infty}h_{1x}^T\big(t,x_*(t)\big) l_1;
           \end{equation}
\item[(b)] for almost all $t \in \R_+$ the maximum condition holds:    
           \begin{equation}\label{PMP3}
           H\big(t,x_*(t),u_*(t),p(t),\lambda_0\big) = \max_{u \in U} H\big(t,x_*(t),u,p(t),\lambda_0\big).
           \end{equation}
\end{enumerate}
\end{theorem}

\begin{example} \label{BeipielRWUnendlich2}
In the Example \ref{BeipielRWUnendlich1} we replace the state constraint $z(t) \leq Z$ by the boundary condition
$$\lim_{t \to \infty} z(t) = Z, \qquad Z> \frac{1}{\varrho}.$$
The bang-bang type control delivers the admissible process $\big(x_*(\cdot),z_*(\cdot),u_*(\cdot)\big)$ again.
Applying the necessary conditions of Theorem \ref{SatzPMP}, we obtain $\lambda_0=1$, $p(t)= e^{-\varrho t}$, $q(t)= \varrho-1$.
Herein, the adjoint $q(\cdot)$ does not vanish at infinity.
\end{example}
\section{The Proof of the Pontryagin Maximum Principle} \label{SectionProof}
\subsection{The Extremal Principle for Locally Convex Problems}
Let $\mathscr{X}$ and $\mathscr{Y}$ be Banach spaces,
let $\mathscr{U}$ be an arbitrary set, let $J$ be a functional on $\mathscr{X} \times \mathscr{U}$,
let $f_1,...,f_l$ be functionals on $\mathscr{X}$,
and let $\mathscr{F}: \mathscr{X} \times \mathscr{U} \to \mathscr{Y}$ be a mapping of the product $\mathscr{X} \times \mathscr{U}$ into $\mathscr{Y}$.
We consider an extremal problem of the abstract form
\begin{equation} \label{Extremalaufgabe}
J(x,u) \to \inf; \quad \mathscr{F}(x,u)=0, \quad f_1(x) \leq 0,...,f_l(x)\leq 0,\quad x \in \mathscr{X},\; u \in \mathscr{U}.
\end{equation}
An admissible pair $(x,u) \in \mathscr{X} \times \mathscr{U}$ in problem (\ref{Extremalaufgabe}) satisfies $\mathscr{F}(x,u)=0$,
$f_1(x) \leq 0,...,f_l(x)\leq 0$ and making the functional $J$ finite.
The admissible pair $(x_*,u_*)$ is a strong local minimizer
iff there exists a number $\varepsilon > 0$ such that the inequality
$J(x,u) \geq J(x_*,u_*)$
holds for any admissible pair $(x,u)$ for which $\|x-x_*\|_{\mathscr{X}} \leq \varepsilon$. \\[2mm]
For the extremal problem (\ref{Extremalaufgabe}) we define the Lagrange function
$\mathscr{L}:\mathscr{X} \times \mathscr{U} \times \R \times \mathscr{Y}^* \times \R^l \to \R$,
$$\mathscr{L}(x,u,\lambda_0,y^*,\lambda_1,...,\lambda_l) = \lambda_0 J(x,u) + \langle y^*, \mathscr{F}(x,u) \rangle + \sum_{j=1}^l \lambda_j f_j(x).$$
Moreover, $\Sigma(\Delta)$ denotes the set
$$\Sigma(\Delta) = \bigg\{ \alpha = (\alpha_1,...,\alpha_{d}) \in \R^{d} \,\bigg|\,
                           \alpha_1,...,\alpha_{d} \geq 0,\, \sum_{i=1}^{d} \alpha_i \leq \Delta \bigg\}.$$
                           
\begin{theorem} \label{TheoremLMR}
Let $(x_*,u_*)$ be an admissible element in the problem (\ref{Extremalaufgabe}).
Assume that the point $x_*$ has a neighborhood $V$ such that
\begin{enumerate}
\item[(a$_1$)] the operator $\mathscr{F}_x(x_*,u_*): \mathscr{X} \to \mathscr{Y}$ is linear and continuous;
\item[(a$_2$)] the function $x \to J(x,u)$ is continuous on the neighborhood $V$ and Fr\'echet differentiable at $x_*$ for any $u \in \mathscr{U}$;
\item[(a$_3$)] the functions $f_1(x),...,f_l(x)$ are continuous on the neighborhood $V$ and regularly locally convex at $x_*$;
\item[(b)] for every finite set of points $u_1,...,u_d \in \mathscr{U}$ and every $\delta >0$,
           there exist a neighborhood $V' \subseteq V$ ($x_* \in V'$), a number $\Delta >0$,
           and a mapping $u:\Sigma(\Delta) \to \mathscr{U}$ which have the following properties
           \begin{enumerate}
           \item[(b$_1$)] $u(0)=u_*$;
           \item[(b$_2$)] the following inequalities hold for any $x,x' \in V'$ and any $\alpha,\alpha' \in \Sigma(\Delta)$:
           \begin{eqnarray*}
           \bigg\| \mathscr{F}\big(x,u(\alpha)\big) - \mathscr{F}\big(x',u(\alpha')\big)-\mathscr{F}_x(x_*,u_*)(x-x') && \\
           -\sum_{k=1}^d (\alpha_k - \alpha_k') \big( \mathscr{F}(x_*,u_k)-\mathscr{F}(x_*,u_*)\big) \bigg\|_{\mathscr{Y}}
              &\leq& \delta \bigg( \|x-x'\|_{\mathscr{X}} + \sum_{k=1}^d |\alpha_k - \alpha_k'|\bigg), \\
           J\big(x,u(\alpha)\big) - J(x,u_*) - \sum_{k=1}^d \alpha_k \big( J(x,u_k)-J(x,u_*)\big)
              &\leq& \delta \bigg( \|x-x_*\|_{\mathscr{X}} + \sum_{k=1}^d \alpha_k \bigg).
                          \end{eqnarray*}
                          
           \end{enumerate}
\end{enumerate}
Finally, we assume that
\begin{enumerate}
\item[(c)] the range of the linear operator $x \to \mathscr{F}_x(x_*,u_*) x$ has a finite codimension in $\mathscr{Y}$.
\end{enumerate}
If $(x_*,u_*)$ is a strong local minimizer in the problem (\ref{Extremalaufgabe}),
then there exist non-trivial Lagrange multipliers
$\lambda_0 \geq 0$, $\lambda_1 \geq 0,...,\lambda_l \geq 0$ and $y^* \in \mathscr{Y}^*$ with
\begin{eqnarray}
&& \label{SatzLMR1} 0 \in \partial_x \mathscr{L}(x_*,u_*,\lambda_0,y^*,\lambda_1,...,\lambda_l)
                    = \lambda_0 J_x(x_*,u_*) + \mathscr{F}_x^*(x_*,u_*)y^*+\sum_{j=1}^l \lambda_j \partial_x f_j(x_*), \\
&& \label{SatzLMR2} \mathscr{L}(x_*,u_*,\lambda_0,y^*,\lambda_1,...,\lambda_l)
                    = \min_{u \in \mathscr{U}} \mathscr{L}(x_*,u,\lambda_0,y^*,\lambda_1,...,\lambda_l), \\
&& \label{SatzLMR3} \lambda_j f_j(x_*)=0, \quad j=1,...,l.
\end{eqnarray}
\end{theorem}

\begin{remark} In view of the extremal problem for regularly locally convex problems of Ioffe \& Tichomirov in \cite{Ioffe},
the extremal principle stated above possesses simplifications:
\begin{enumerate}
\item[(1)] The functional $x \to J(x,u)$ assumed to be Fr\'echet differentiable at $x_*$,
           which is a stronger requirement than the regularly locally convexity in \cite{Ioffe}.
\item[(2)] In \cite{Ioffe}, the mapping in assumption (b) has the more abstract form
           $u(x,\alpha):V' \times \Sigma(\Delta) \to U$.
\end{enumerate}
\end{remark}

\begin{remark} Since the first inequality in the assumption (b$_2$) states the required differentiability of the mapping $\mathscr{F}$,
we replaced the assumption (a) in \cite{Ioffe}, e.g.
the mapping $x \to \mathscr{F}(x,u)$ is of class $C_1$ at the point $x_*$ for any $u \in \mathscr{U}$,
by the assumption (a$_1$).
The type of differentiability in (b$_2$) satisfies the requirements of the implicit function theorem in form of the generalized Lyusternik theorem in \cite{Ioffe}.
\end{remark}
\subsection{The Space of Continuous Functions Converging at Infinity} \label{SectionSpace}
In the following we develop a framework for the state space $\mathscr{X}$ in the abstract extremal problem (\ref{Extremalaufgabe}).
The need of a suitable state space arises by different challenges on the infinite horizon:
\begin{enumerate}
\item According to (\ref{Aufgabe3}), e.\,g. $\lim\limits_{t \to \infty} h_1\big(t,x(t)\big)=0$,
      the trajectory $x(\cdot)$ may posses a limit at infinity.
\item Let $C_0(\R_+,\R)$ be the space of continuous functions vanishing at infinity.
      Then $x(\cdot) \in C_0(\R_+,\R)$ satisfies the state constraint $x(t) \leq 0$ on $\R_+$ iff
      $$x(\cdot) \in \mathscr{K}=\{ z(\cdot) \in C_0(\R_+,\R) \,|\, z(t)\leq 0 \mbox{ for all } t \in \R_+\}.$$
      But $\mathscr{K}$ possesses an empty interior:
      Let $z(\cdot) \in \mathscr{K}$ with $z(t) <0$ on $\R_+$ and consider the sequence
      $$z_n(t)=\left\{\begin{array}{ll} z(t), & t \not\in [2n-1,2n+1], \\
                                        z(t)+2|z(2n)| \cdot (t-(2n-1))((2n+1)-t), & t \in [2n-1,2n+1].
                      \end{array} \right.$$
      Then $\|z_n(\cdot)-z(\cdot)\|_\infty \to 0$ as $n \to \infty$ and $z_n(\cdot) \not\in \mathscr{K}$ for all $n$.
      Hence in the space $C_0(\R_+,\R)$ the interior of $\mathscr{K}$ and the origin cannot be separated by the
      Hahn-Banach separation theorem for convex sets. 
\item In the space $C_0(\R_+,\R)$ the subdifferential of
      $f\big(x(\cdot)\big)=\sup\limits_{t \in \R_+} x(t)=\max\limits_{t \in \overline{\R}_+} x(t)$
      delivers in $x(\cdot)=0$ the set of all non-negative regular Borel measures $\mu$ on $\R_+$ with the
      non-strict attribute $\|\mu\| \leq 1$.
      Indeed, for any function $z(\cdot) \in C_0(\R_+,\R)$ with $z(t) <0$ on $\R_+$ the definition 
      $$\partial f\big(z(\cdot)\big) =\{x^* \in \partial f(0) \,|\, f\big(z(\cdot)\big)= \langle x^*,z(\cdot) \rangle \}$$
      leads to
      $$0 = \sup_{t \in \R_+} z(t) =\int_0^\infty z(t) \, d\mu(t) \quad\Leftrightarrow\quad \mu=0.$$
\end{enumerate}

To fix these lacks we consider the space $C_{\lim}(\R_+,\R^n)$ of continuous vector-functions converging at infinity,
which can be defined equivalently in the following ways:
\begin{eqnarray*}
C_{\lim}(\R_+,\R^n) &=& \{x(\cdot) \in C(\R_+,\R^n) \,|\, \lim_{t \to \infty} x(t)= a \in \R^n\} \\
                    &=& \{x(\cdot): \R_+ \to \R^n \,|\, x(t)=x_0(t)+a, x_0(\cdot) \in C_0(\R_+,\R^n), a \in \R^n\}.
\end{eqnarray*}

The first setting delivers the ''natural'' norm $\|x(\cdot)\|=\|x(\cdot)\|_\infty$.
But the representation in the second setting implies the norm $\|x(\cdot)\|=\|x_0(\cdot)\|_\infty + \|a\|$.

\begin{lemma} \label{LemmaNorm}
On $X=C_{\lim}(\R_+,\R^n)$ the norms $\|x(\cdot)\|_X=\|x(\cdot)\|_\infty$ and $\|x(\cdot)\|_X=\|x_0(\cdot)\|_\infty + \|a\|$
are equivalent.
\end{lemma}

\textbf{Proof.} By definition,
any $x(\cdot) \in C_{\lim}(\R_+,\R^n)$ possesses the unique representation $x(\cdot)=x_0(\cdot)+a$
with $x_0(\cdot) \in C_0(\R_+,\R^n)$ and $a \in \R^n$.
Since $x(t) \to a$ as $t \to \infty$, the inequality $\|x(\cdot)\|_\infty \geq \|a\|$ holds.
Let $\|x_0(\cdot)\|_\infty \leq \|a\|$.
Then we obtain immediately $\|x(\cdot)\|_X \leq 2\|x(\cdot)\|_\infty$.
In the case $\|a\| \leq \|x_0(\cdot)\|_\infty$ there exists $\lambda \in [0,1]$ with $\|a\| = \lambda \|x_0(\cdot)\|_\infty$.
It follows
$$\|x(\cdot)\|_\infty = \|x_0(\cdot)+a\|_\infty \geq \max \{\|x_0(\cdot)\|_\infty - \|a\|,\|a\|\}
  \geq \min_{\lambda \in [0,1]} \max\{1-\lambda,\lambda\}\cdot \|x_0(\cdot)\|_\infty = \frac{1}{2} \|x_0(\cdot)\|_\infty.$$
Consequently, the inequalities $\|x(\cdot)\|_\infty \leq \|x(\cdot)\|_X \leq 3\|x(\cdot)\|_\infty$ hold on
$C_{\lim}(\R_+,\R^n)$. \hfill $\square$ \\

In the following we make use of the notation $x(\infty)=\lim\limits_{t \to \infty} x(t)$
instead of the vector $a \in \R^n$.

\begin{lemma}[Riesz's representation theorem] \label{CorollaryDarstellung}
Any continuous linear functional $x^*$ on $C_{\lim}(\R_+,\R^n)$ can be uniquely represented in the form
$$\langle x^*(\cdot),x(\cdot) \rangle = \int_0^\infty \langle x(t)-x(\infty), d\nu(t) \rangle
              +\alpha^T x(\infty) = \int_0^\infty \langle x_0(t), d\nu(t) \rangle +\alpha^T x(\infty),$$
where $\nu=(\nu_1,...,\nu_n)$ is a vector of regular signed Borel measures on $\R_+$ (cf. \cite{Rudin}) and $\alpha \in \R^n$,
or equivalently, in the form
$$\langle x^*(\cdot),x(\cdot) \rangle = \int_0^\infty \langle x(t), d\mu(t) \rangle,$$
where $\mu=(\nu_1,...,\nu_n)$ is a vector of signed Borel measures $\mu_1,...,\mu_n \in \mathscr{M}(\overline{\R}_+)$
(cf. Definition \ref{DefinitionBorel}).
\end{lemma}

\textbf{Proof.} Any $x(\cdot) \in C_{\lim}(\R_+,\R^n)$ has the unique representation $x(\cdot)=x_0(\cdot)+a$
with $x_0(\cdot) \in C_0(\R_+,\R^n)$ and $a \in \R^n$.
Consider the linear mapping $\Lambda$ with $\Lambda x(\cdot) = x_0(\cdot)$,
which maps the space $C_{\lim}(\R_+,\R^n)$ onto $C_0(\R_+,\R^n)$.
Then $\Lambda$ is continuous ($\|x_0(\cdot)\|_\infty \leq 3\|x(\cdot)\|_\infty$ by Lemma \ref{LemmaNorm}) and
$$Ker\,\Lambda=\{ x(\cdot) \in C_{\lim}(\R_+,\R^n \,|\, x(t)=a \mbox{ on } \R_+\}.$$
Let $x^* \in C^*_{\lim}(\R_+,\R^n)$.
We denote by $\alpha_i$ the value of the functional $x^*$ at the vector-valued function whose $i$-th component is identical $1$ and
whose remaining components are identically zero.
Now, we consider the functional $x_1^* \in C^*_{\lim}(\R_+,\R^n)$ defined by the formula
$$\langle x_1^*,x(\cdot) \rangle = \langle x^*,x(\cdot) \rangle - \alpha^T a, \quad x(t)=x_0(t)+a, \quad \alpha=(\alpha_1,...,\alpha_n).$$
Obviously, $x_1^* \in (Ker\,\Lambda)^\perp$.
By the lemma of the annihilator (cf. \cite{Ioffe}),
there exists a functional $x_0^* \in X_0^*$ with $x_1^* = \Lambda^* x_0^*$.
That means the equation
$$\langle x_1^*,x(\cdot) \rangle = \langle \Lambda^* x_0^*,x(\cdot) \rangle = \langle x_0^*,\Lambda x(\cdot)\rangle =\langle x_0^*,x_0(\cdot) \rangle$$
holds for all $x(\cdot) \in C_{\lim}(\R_+,\R^n)$.
Consequently, for $x^* \in C^*_{\lim}(\R_+,\R^n)$ we obtain the representation
$$\langle x^*,x(\cdot) \rangle =  \langle x_0^*, x_0(\cdot) \rangle + \alpha^T a
   = \int_0^\infty \langle x_0(t), d\nu(t) \rangle +\alpha^T x(\infty).$$
In this equation we set
$\mu_\infty=\mu(\{\infty\})=\alpha -\displaystyle\int_0^\infty d\nu(t)$ and define the Borel measure $\mu=\nu+\mu_\infty$.
Then, the last representation formula in Lemma \ref{CorollaryDarstellung} is shown.
Finally, the uniqueness of this representation can be verified directly. \hfill $\square$ \\

In the space $C_{\lim}(\R_+,\R^n)$ we conclude:
\begin{enumerate}
\item The space $C_{\lim}(\R_+,\R^n)$ is a suitable framework for restrictions of the form
      $\lim\limits_{t \to \infty} h_1\big(t,x(t)\big)=0$.
\item The interior of the cone $\mathscr{K}=\{ z(\cdot) \in C_{\lim}(\R_+,\R) \,|\, z(t)\leq 0 \mbox{ for all } t \in \R_+\}$
      is non-empty and containing any $x(\cdot) \in C_{\lim}(\R_+,\R)$ with
      $\max\limits_{t \in \overline{\R}_+} x(t) <0$.      
\item We consider $f\big(x(\cdot)\big)=\max\limits_{t \in \overline{\R}_+} x(t)$.
      By definition of $\partial f(0)$ we obtain in the space $C_{\lim}(\R_+,\R)$ the inequalities
      $$\max_{t \in \overline{\R}_+} x(t) \geq \int_0^\infty x(t) \, d\mu(t)
        \geq - \max_{t \in \overline{\R}_+} \big( -x(t)\big)
        = \min_{t \in \overline{\R}_+} x(t).$$ 
      This shows, that the subdifferential of $f$ consists in $x(\cdot)=0$ of those and only those non-negative Borel measures
      $\mu \in \mathscr{M}(\overline{\R}_+)$, which satisfy $\|\mu\| = 1$.
\end{enumerate}
\subsection{The Abstract Optimization Problem and Properties} \label{SectionEP}
Let $\big(x_*(\cdot),u_*(\cdot)\big) \in \mathscr{A}_{Lip} \cap \mathscr{A}_{lim}$ and $x_*(\cdot) \in C_{\lim}(\R_+,\R^n)$.
On $C_{\lim}(\R_+,\R^n) \times L_\infty(\R_+,U)$ we consider the mappings
\begin{eqnarray*}
&& J\big(x(\cdot),u(\cdot)\big) = \int_0^\infty \omega(t) f\big(t,x(t),u(t)\big) \, dt, \\
&& F\big(x(\cdot),u(\cdot)\big)(t) = x(t)-x(0) -\int_0^t \varphi\big(s,x(s),u(s)\big) \, ds, \quad t \geq 0, \\
&& H_0\big(x(\cdot)\big)= h_0\big(x(0)\big), \qquad H_1\big(x(\cdot)\big)= \lim_{t \to \infty} h_1\big(t,x(t)\big), \\
&& G_j\big(x(\cdot)\big)(t)= g_j\big(t,x(t)\big), \; t \in \overline{\R}_+,
   \qquad f_j\big(x(\cdot)\big) = \max_{t \in \overline{\R}_+} g_j\big(t,x(t)\big), \qquad j=1,...,l.
\end{eqnarray*}
Moreover, on $C_{\lim}(\R_+,\R^n)$ we define the linear operators
\begin{eqnarray*}
&& J_x\big(x(\cdot),u(\cdot)\big) \xi(\cdot)
   = \int_0^\infty \omega(t) \langle f_x\big(t,x(t),u(t)\big) , \xi(t) \rangle \, dt, \\
&& \big[F_x\big(x(\cdot),u(\cdot)\big) \xi(\cdot) \big](t)
   = \xi(t)-\xi(0) -\int_0^t \varphi_x\big(s,x(s),u(s)\big)\xi(s) \, ds, \quad t \geq 0, \\
&& H_0'\big(x(\cdot)\big) \xi(\cdot)= \langle h_0'\big(x(0)\big) , \xi(0) \rangle,
   \qquad H_1'\big(x(\cdot)\big) \xi(\cdot)= \lim_{t \to \infty} \langle h_{1x}\big(t,x(t)\big) , \xi(t) \rangle, \\
&& \big[G_j'\big(x_*(\cdot)\big)x(\cdot)\big](t)= \langle g_{jx}\big(t,x_*(t)\big), x(t) \rangle, \quad t \in \overline{\R}_+, \quad j=1,...,l.
\end{eqnarray*}
Furthermore, let $\mathscr{X}=C_{\lim}(\R_+,\R^n)$ and let $\mathscr{U}$ be the set
\begin{eqnarray}
&& \mathscr{U} = \{ u(\cdot) \in L_\infty(\R_+,U) \,|\, u(t)=u_*(t) + \chi_{M}(t)\big(w(t)-u_*(t)\big),
                    w(\cdot) \in L_\infty(\R_+,U), \nonumber\\
\label{DefinitionU} && \hspace*{35mm} M \subset \R_+ \mbox{ measurable and bounded} \}.
\end{eqnarray}
We introduce the mappings
$$\mathscr{F}=(F,H_0,H_1), \quad \mathscr{F}_x\big(x_*(\cdot),u_*(\cdot)\big)=(F_x,H_0',H_1')\big(x_*(\cdot),u_*(\cdot)\big).$$
Then the optimal control problem (\ref{Aufgabe1})--(\ref{Aufgabe5}) becomes the form (\ref{Extremalaufgabe}):
\begin{equation} \label{Extremalaufgabe1}
J\big(x(\cdot),u(\cdot)\big) \to \inf; \quad \mathscr{F}\big(x(\cdot),u(\cdot)\big)=0,\quad
f_1\big(x(\cdot)\big) \leq 0,...,f_l\big(x(\cdot)\big)\leq 0, \quad 
x(\cdot) \in \mathscr{X},\; u(\cdot) \in \mathscr{U}.
\end{equation}
\begin{lemma} \label{LemmaZielfunktional}
For any $u(\cdot) \in \mathscr{U}$ the mapping $x(\cdot) \to J\big(x(\cdot),u(\cdot)\big)$ is Fr\'echet differentiable at the
point $x_*(\cdot)$, and $J_x\big(x_*(\cdot),u(\cdot)\big)$ is the Fr\'echet derivative.
\end{lemma}

\textbf{Proof.} According to $x_*(\cdot) \in \mathscr{X}_{Lip}$,
for a fixed $u(\cdot) \in L_\infty(\R_+,U)$
the mappings $t \to f\big(t,x(t),u(t)\big)$ and $t \to f_x\big(t,x(t),u(t)\big)$ are uniformly bounded
by $C_0>0$ a.e. on $\R_+$ for all $x(\cdot) \in C_{\lim}(\R_+,\R^n)$ with $\|x(\cdot)-x_*(\cdot)\|_\infty \leq \gamma$.
Therefore, the functional $J\big(x(\cdot),u(\cdot)\big)$ is finite.
Moreover, the mapping $x(\cdot) \to J_x\big(x_*(\cdot),u(\cdot)\big) x(\cdot)$ is linear and continuous.
Let $\varepsilon >0$ be given.
Then we choose a number $T>0$ with
$$\int_T^\infty \omega(t) 2 C_0 \, dt \leq \frac{\varepsilon}{2}.$$
Moreover, for sufficiently small $\lambda_0 >0$ we have
$$\esssup_{t \in [0,T], \|x\|_\infty \leq 1} \big\|f_x\big(t,x_*(t)+\lambda x,u(t)\big)-f_x\big(t,x_*(t),u(t)\big)\big\|
  \leq \frac{\varepsilon}{2} \frac{1}{\|\omega(\cdot)\|_{L_1}}, \quad 0<\lambda\leq \lambda_0.$$
It follows
\begin{eqnarray*}
&& \bigg|\frac{J\big(x_*(\cdot)+\lambda x(\cdot),u(\cdot)\big) - J\big(x_*(\cdot),u(\cdot)\big)}{\lambda} - J_x\big(x_*(\cdot),u(\cdot)\big) x(\cdot)\bigg| \\
&=& \bigg|\int_0^\infty \omega(t)
         \bigg[ \int_0^1 \langle f_x\big(t,x_*(t)+\lambda s x(t),u(t)\big) - f_x\big(t,x_*(t),u(t)\big), x(t) \rangle ds \bigg] \, dt \bigg| \\
&\leq& \int_0^T \omega(t) \frac{\varepsilon}{2\|\omega(\cdot)\|_{L_1}} \|x(t)\| \, dt + \int_T^\infty \omega(t) 2 C_0 \|x(t)\| dt \leq \varepsilon
\end{eqnarray*}
for all $x(\cdot) \in C_{\lim}(\R_+,\R^n)$ with $\|x(\cdot)\|_\infty \leq 1$ and all
$0 < \lambda \leq \lambda_0$. \hfill $\square$ \\

\begin{lemma} \label{LemmaDynamik1}
For any $u(\cdot) \in \mathscr{U}$,
the operator $x(\cdot) \to F\big(x(\cdot),u(\cdot)\big)$ maps a neighborhood of $x_*(\cdot)$ into the space $C_{\lim}(\R_+,\R^n)$.
\end{lemma}

\textbf{Proof.} According to (\ref{condition1}), (\ref{condition2}) and $u(\cdot) \in \mathscr{U}$,
there exist $\gamma, T>0$ such that $u(t)=u_*(t)$ for all $t \geq T$ and
\begin{eqnarray*}
       \int_T^\infty \big\| \varphi\big(t,x(t),u_*(t)\big)\big\| \, dt 
&\leq& \|x(\cdot)-x_*(\cdot)\|_\infty + \int_T^\infty \big\| \varphi\big(t,x_*(t),u_*(t)\big) \big\| \, dt \\
&    & + \int_T^\infty \big\| \varphi_x\big(t,x_*(t),u_*(t)\big) \big(x(t)-x_*(t)\big) \big\| \, dt < \infty
\end{eqnarray*}
for all $x(\cdot) \in C_{\lim}(\R_+,\R^n)$ with $\|x(\cdot)-x_*(\cdot)\|_\infty \leq \gamma$.
Furthermore, the condition
$$\int_0^T \big\| \varphi\big(t,x(t),u(t)\big)\big\| \, dt< \infty$$
holds.
Consequently, $\lim\limits_{t \to \infty} F\big(x(\cdot),u(\cdot)\big)(t)$ exists
for any pair $\big(x(\cdot),u(\cdot)\big) \in C_{\lim}(\R_+,\R^n) \times \mathscr{U}$. \hfill $\square$

\begin{lemma} \label{LemmaDynamik2}
The linear operator $\mathscr{F}_x\big(x_*(\cdot),u_*(\cdot)\big)$ is continuous.
\end{lemma}

\textbf{Proof.} At first it is obvious that the mappings $H_i:C_{\lim}(\R_+,\R^n) \to \R^{s_i}$ are continuously 
differentiable on $C_{\lim}(\R_+,\R^n)$
and the linear operators $H_i'\big(x(\cdot)\big)$ are their derivatives.
According to the second condition in (\ref{condition1}),
the linear operator $F_x\big(x_*(\cdot),u_*(\cdot)\big)$ is continuous. \hfill $\square$ \\

\begin{lemma} \label{LemmaZB1}
For $j=1,...,l$ the mappings $G_j$ are differentiable at the point $x_*(\cdot)$ with Fr\'echet derivative
$G_j'\big(x_*(\cdot)\big)$.
\end{lemma}

\textbf{Proof.} According to the properties of $g_j(t,x)$,
the limits
$\lim\limits_{t \to \infty} g_j\big(t,x(t)\big)$, $\lim\limits_{t \to \infty} g_{jx}\big(t,x(t)\big)$
exist for all $x(\cdot) \in C_{\lim}(\R_+,\R^n)$ with $\|x(\cdot)-x_*(\cdot)\|_\infty \leq \gamma$. 
Moreover, any $G_j'\big(x_*(\cdot)\big)$ is linear and continuous.
Let $\varepsilon >0$ be given.
Then, for sufficiently small $\lambda_0 >0$, the inequality
$$\max_{t \in \overline{\R}_+, \|x\|_\infty \leq 1} \big\|g_{jx}\big(t,x_*(t)+\lambda x\big)-g_{jx}\big(t,x_*(t)\big)\big\| \leq \varepsilon$$
holds for all $\lambda$ with $0 < \lambda \leq \lambda_0$.
Therefore, for any $j=1,...,l$ the relation
\begin{eqnarray*}
\lefteqn{\bigg\|\frac{G_j\big(x_*(\cdot)+\lambda x(\cdot)\big) - G_j\big(x_*(\cdot)\big)}{\lambda}
               - G_j'\big(x_*(\cdot)\big) x(\cdot)\bigg\|_\infty} \\
&&= \max_{t \in \overline{\R}_+}
    \bigg|\int_0^1 \langle g_{jx}\big(t,x_*(t)+\lambda s x\big)-g_{jx}\big(t,x_*(t)\big), x(t) \rangle ds \bigg|
    \leq \varepsilon
\end{eqnarray*}
is satisfied for all $x(\cdot) \in C_{\lim}(\R_+,\R^n)$ with $\|x(\cdot)-x_*(\cdot)\|_\infty \leq 1$
and all $0 < \lambda \leq \lambda_0$. \hfill $\square$
\subsection{Generalized Needle Variations}
Introducing needle variations in infinite horizon optimal control we have to take into account two challenging tasks:
the unbounded time interval and an unbounded density $\omega(\cdot) \in L_1(\R_+,\R)$
like a Weibull type distribution $\omega(t)=t^{k-1}e^{-t^k}$ with form parameter $k \in (0,1)$.
In the following we develop a method,
where the needle variations are concentrated on a suitable compact subset of $\R_+$. \\
For any measurable subset $A$ the number $|A|$ denotes the Lebesgue measure of $A$.
Let $K \subset \R_+$ be a compact set with $|K|>0$.
Then there exists a number $T>0$ with $K \subseteq [0,T]$.
Moreover, for any $n \in \N$ there exist numbers $0=t_0 < t_1 <... < t_n =T$ with
$|[t_{i-1},t_i] \cap K| = |K| / n$, $i=1,...,n$.
Therefore, the sets $\Delta_i= [t_{i-1},t_i] \cap K$ defining a sorted partition of $K$ with equal length.
With these preliminary remarks we can state the following Lemma,
which is a direct consequence of the results on intervals $[t_0,t_1]$ in \cite{Ioffe}, pp. 243--246:

\begin{lemma} \label{LemmaNV}
Let $K$ be a compact subset of $\R_+$ and let $y_i(\cdot),\, y_i : K \to \R^{n_i},\, i=1,...,d,$ be bounded measurable vector-valued functions.
Then, for every $\delta > 0$, there exist one-parameter families $M_1(\alpha),...,M_d(\alpha)$, $0\leq \alpha \leq 1/d$,
of measurable subsets of $K$,
such that
\begin{eqnarray*}
&& |M_i(\alpha)| = \alpha |K|, \qquad M_i(\alpha') \subseteq M_i(\alpha), \qquad M_i(\alpha) \cap M_{i'}(\alpha') = \emptyset,\\
&& \max_{t \in K}
   \bigg\| \int_{[0,t] \cap K} \big( \chi_{M_i(\alpha)}(\tau) - \chi_{M_i(\alpha')}(\tau) \big) y_i(\tau) \, d\tau 
            - (\alpha-\alpha') \int_{[0,t] \cap K} y_i(\tau) \, d\tau \bigg\| \leq \delta |\alpha-\alpha'|
\end{eqnarray*}
for all $i=1,...,d$, $0 \leq \alpha' \leq \alpha \leq 1/d$ and $i \not= i'$.
\end{lemma}

\textbf{Proof.} The proof in \cite{Ioffe} bases on the explicit construction for step functions.
We repeat this construction in the case of the compact set $K$.
Let $y(t)=\sum\limits_{j=1}^m y_j \chi_{A_j}(t)$ be a step function,
where $\{A_1,...,A_m\}$ is a partition of $K$.
Furthermore, let $C= \max\limits_j \|y_j\|$.
Then we can state a sorted partition $\{\Delta_1,...,\Delta_r\}$ of $K$ with equal length no greater than $\delta / (2C)$. \\
We define the sets
$$M_{ij}(\alpha) = \bigg\{ t \in (A_j \cap \Delta_i)
                      \,\bigg|\, \int_0^t \chi_{(A_j \cap \Delta_i)} (\tau) \, d\tau
                                                                < \alpha \cdot |A_j \cap \Delta_i| \bigg\}.$$
Then the sets $M(\alpha) = \displaystyle\bigcup_{i,j} M_{ij}(\alpha)$ are the required sets.
It holds
$$|M(\alpha)| = \sum_{i,j} |M_{ij}(\alpha)| = \alpha \sum_{i,j} |A_j \cap \Delta_i| = \alpha |K|.$$
If $\alpha \geq \alpha'$, then $M_{ij}(\alpha') \subseteq M_{ij}(\alpha)$ and, therefore, $M(\alpha') \subseteq M(\alpha)$.
Finally,
$$\int_{M_{ij}(\alpha)} y(t) \, dt = \int_K \chi_{M_{ij}(\alpha)}(t) y(t) \, dt = \alpha y_j \cdot |A_j \cap \Delta_i|.$$
It follows that the values of $\displaystyle \int_{[0,t] \cap K} (\alpha-\alpha') y(\tau) \, d\tau$ and
$\displaystyle \int_{[0,t] \cap K} \big( \chi_{M(\alpha)}(\tau) - \chi_{M(\alpha')}(\tau) \big) y(\tau) \, d\tau$
coincide at the end of the subintervals $\Delta_i$:
\begin{eqnarray*}
\lefteqn{\int_{\Delta_i} \big( \chi_{M(\alpha)}(t) - \chi_{M(\alpha')}(t) \big) y(t) \, dt
         = \sum_j \left( \int_{M_{ij}(\alpha)} y(t) \, dt - \int_{M_{ij}(\alpha')} y(t) \, dt \right) } \\
&=& (\alpha-\alpha') \sum_j y_j \cdot |A_j \cap \Delta_i| 
    = (\alpha-\alpha') \sum_j \int_{\Delta_i} \chi_{(A_j \cap \Delta_i)}(t) y(t) \, dt
    = (\alpha-\alpha') \int_{\Delta_i} y(t) \, dt.
\end{eqnarray*}
Moreover, if $t \in \Delta_i$, then
\begin{eqnarray*}
\lefteqn{\left\| \int_{[0,t] \cap \Delta_i} \Big[ \big( \chi_{M(\alpha)}(\tau) - \chi_{M(\alpha')}(\tau) \big) y(\tau)
                                 - (\alpha-\alpha') y(\tau) \Big] \, d\tau \right\|} \\
&\leq& \left\| \int_{[0,t] \cap \Delta_i}
               \big( \chi_{M(\alpha)}(\tau) - \chi_{M(\alpha')}(\tau) \big) y(\tau) \, d\tau \right\|
       + |\alpha-\alpha'| \left\| \int_{[0,t] \cap \Delta_i} y(\tau) \, d\tau \right\| \\
&\leq& 2 C |\alpha-\alpha'| \cdot |\Delta_i| \leq \delta |\alpha-\alpha'|.
\end{eqnarray*}
The construction for step functions is completed. \hfill $\square$ \\

We return to the problem (\ref{Extremalaufgabe1}).
We have to verify that this problem has the properties indicated in the condition b) of Theorem \ref{TheoremLMR}.
Let $\big(x_*(\cdot),u_*(\cdot)\big) \in \mathscr{A}_{adm} \cap \mathscr{A}_{Lip} \cap \mathscr{A}_{lim}$.
Moreover,
let $u_1(\cdot),...,u_d(\cdot) \in \mathscr{U}$ and $\delta >0$ be given.
Then there exists a number $T>0$,
such that the bounded sets $M_i$ in the representation of the elements $u_i(\cdot)$ are subsets of the interval $[0,T]$.
In addition, the number $T$ can be choosen,
such that the condition (\ref{condition2}) is satisfied with $\delta/3$, e.g.
\begin{equation} \label{NVCondition1} 
\int_T^\infty \big\| \varphi\big(t,x(t),u_*(t)\big)-\varphi\big(t,x'(t),u_*(t)\big) - \varphi_x\big(t,x_*(t),u_*(t)\big)\big(x(t)-x'(t)\big) \big\| \, dt
\leq \frac{\delta}{3} \|x(\cdot)-x'(\cdot)\|_\infty
\end{equation}
holds for all $x(\cdot), x'(\cdot) \in W^1_\infty(\R_+,\R^n)$ with $\|x(\cdot)-x_*(\cdot)\|_\infty \leq \gamma$,
$\|x'(\cdot)-x_*(\cdot)\|_\infty \leq \gamma$.
Since $\big(x_*(\cdot),u_*(\cdot)\big) \in \mathscr{A}_{Lip}$,
there exists $0<\sigma \leq \gamma$, such that
\begin{equation} \label{NVCondition2} 
\int_0^T \big\| \varphi\big(t,x(t),u_*(t)\big)-\varphi\big(t,x'(t),u_*(t)\big) - \varphi_x\big(t,x_*(t),u_*(t)\big)\big(x(t)-x'(t)\big) \big\| \, dt
\leq \frac{\delta}{3} \|x(\cdot)-x'(\cdot)\|_\infty
\end{equation}
holds for all $x(\cdot), x'(\cdot) \in W^1_\infty(\R_+,\R^n)$ with $\|x(\cdot)-x_*(\cdot)\|_\infty \leq \sigma$,
$\|x'(\cdot)-x_*(\cdot)\|_\infty \leq \sigma$.
Furthermore, due to Lusin's theorem there exists a compact set $K \subseteq [0,T]$,
such that $t \to \omega(t)$ is continuous on $K$ and the relations
\begin{eqnarray}
&& \label{NVCondition3} \int_{[0,T] \setminus K} \big\|\varphi\big(t,x_*(t),u_i(t)\big)-\varphi\big(t,x_*(t),u_*(t)\big)\big\| \, dt \leq \frac{\delta}{3}, \\
&& \label{NVCondition4} \int_{[0,T] \setminus K} \omega(t)\big|f\big(t,x(t),u_i(t)\big)-f\big(t,x(t),u_*(t)\big)\big| \, dt \leq \frac{\delta}{2}
\end{eqnarray}
hold for $i=1,...,d$ and for all $x(\cdot) \in W^1_\infty(\R_+,\R^n)$ with $\|x(\cdot)-x_*(\cdot)\|_\infty < \sigma$. \\
In the following we consider the vector-valued functions $y_i(\cdot)$,
$$y_i (t) = \Big( \varphi\big(t,x_*(t),u_i(t) \big) - \varphi\big(t,x_*(t),u_*(t) \big),
                  \omega(t)\big[ f\big(t,x_*(t),u_i(t) \big) - f\big(t,x_*(t),u_*(t) \big)\big] \Big), \quad i = 1,...,d.$$
On $K$ the mappings $t \to y_i(t)$ are measurable and bounded.
With the corresponding sets $M_1(\alpha) \subseteq K$,...,$M_d(\alpha) \subseteq K$ in Lemma \ref{LemmaNV} we define on $\R_+$ 
the mapping $\alpha \to u_\alpha(\cdot) \in \mathscr{U}$, where
$$u_\alpha(t) = u_*(t) + \sum_{i=1}^{d} \chi_{M_i(\alpha_i)}(t) \cdot \big( u_i(t)-u_*(t) \big).$$
According to this expression,
we call the mapping $\alpha \to u_\alpha(\cdot)$ a generalized needle variation.
On the set $Q^{d}$,
$$Q^{d} = \Big\{ \alpha = (\alpha_1,...,\alpha_{d}) \in \R^{d} \,\Big|\,
                                0 \leq \alpha_i  \leq 1/d, \, i =1,...,d \Big\},$$
the mapping $\alpha \to u_\alpha(\cdot) \in \mathscr{U}$ is well-defined.
By $V(\sigma)$ and $\Sigma(\Delta)$ we denote the sets
\begin{eqnarray*}
V(\sigma) &=& \{ x(\cdot) \in C_{\lim}(\R_+,\R^n) \,\big|\, \|x(\cdot) - x_*(\cdot)\|_{\infty} \leq \sigma \}, \\
\Sigma(\Delta) &=& \left\{ \alpha = (\alpha_1,...,\alpha_{d}) \in \R^{d} \,\bigg|\,
                                \alpha_1,...,\alpha_{d} \geq 0,\, \sum_{i=1}^{d} \alpha_i \leq \Delta \right\}.
\end{eqnarray*}
On $V(\sigma) \times \Sigma(\Delta)$ we consider the mappings
\begin{eqnarray*}
    \Phi_1\big(x(\cdot),\alpha\big)(t)
&=& \int_0^t \big[\varphi\big(\tau,x(\tau),u_\alpha(\tau)\big) - \varphi\big(\tau,x_*(\tau),u_*(\tau)\big)\big] \, d\tau, \quad t \in \R_+, \\
    \Lambda_1\big(x(\cdot),\alpha\big)(t)
&=& \int_0^t \bigg[\varphi_x\big(\tau,x_*(\tau),u_*(\tau)\big) \big(x(\tau)-x_*(\tau)\big) \\
& & \hspace*{7mm} + \sum_{i=1}^{d} \alpha_i \cdot 
    \Big( \varphi\big(\tau,x_*(\tau),u_i(\tau)\big) - \varphi\big(\tau,x_*(\tau),u_*(\tau)\big)\Big)\bigg] \, d\tau, \quad t \in \R_+
\end{eqnarray*}
and the functionals
\begin{eqnarray*}
    \Phi_2\big(x(\cdot),\alpha\big)
&=& \int_0^\infty \omega(t)\big[f\big(t,x(t),u_\alpha(t)\big)-f\big(t,x(t),u_*(t)\big)\big] \, dt, \\
    \Lambda_2\big(x(\cdot),\alpha\big)
&=& \sum_{i=1}^{d} \int_0^\infty \alpha_i \cdot \omega(t) \big[ f\big(t,x(t),u_i(t)\big)
                                                    - f\big(t,x(t),u_*(t)\big)\big] \, dt.
\end{eqnarray*}

\begin{lemma} \label{LemmaNVUH1}
There exist $\sigma>0$ and $\Delta \in (0,1/d]$,
such that the condition
\begin{equation} \label{NVUH1}
\Big\| \big[\Phi_1\big(x(\cdot),\alpha\big)-\Phi_1\big(x'(\cdot),\alpha'\big)
          -\Lambda_1\big(x(\cdot),\alpha\big)+\Lambda_1\big(x'(\cdot),\alpha'\big)\big] (\cdot) \Big\|_\infty
\leq\delta \bigg( \|x(\cdot)-x'(\cdot)\|_\infty + \sum_{i=1}^d |\alpha_i - \alpha_i'| \bigg)
\end{equation}
holds for all $x(\cdot), x'(\cdot) \in V(\sigma)$ and all $\alpha, \alpha' \in \Sigma(\Delta)$.  
\end{lemma}

\textbf{Proof.} We take into account that $u_\alpha(t)=u_*(t)$ for $t \not\in K$ and that $u_i(t)=u_*(t)$ for $t>T$.
Hence, the left-hand side in (\ref{NVUH1}) can be estimated as follows:
\begin{eqnarray*}
& & \int_{[0,T] \setminus K} \big\| \varphi\big(t,x(t),u_*(t)\big)-\varphi\big(t,x'(t),u_*(t)\big)
                         - \varphi_x\big(t,x_*(t),u_*(t)\big)\big(x(t)-x'(t)\big) \big\| \, dt  \\
&+& \sum_{i=1}^d (\alpha_i- \alpha'_i) \cdot \int_{[0,T] \setminus K} \big\|\varphi\big(t,x_*(t),u_i(t)\big)-\varphi\big(t,x_*(t),u_*(t)\big)\big\| \, dt \\
&+& \max_{t \in [0,T]} \bigg\| \int_{[0,t] \cap K} \bigg[\varphi \big( \tau, x(\tau), u_\alpha(\tau) \big)
                            - \varphi \big( \tau, x'(\tau), u_{\alpha'}(\tau) \big) \\
& & \hspace*{30mm} - \varphi_x \big( \tau, x_*(\tau), u_*(\tau) \big) \big( x(\tau) - x'(\tau) \big) \\
& & \hspace*{30mm} - \sum_{i=1}^d (\alpha_i- \alpha'_i) \cdot 
    \Big( \varphi \big( \tau, x_*(\tau), u_i(\tau) \big) - \varphi \big( \tau, x_*(\tau), u_*(\tau) \big) \Big) \bigg] d\tau \bigg\| \\
&+& \int_T^\infty \big\| \varphi\big(t,x(t),u_*(t)\big)-\varphi\big(t,x'(t),u_*(t)\big)
                         - \varphi_x\big(t,x_*(t),u_*(t)\big)\big(x(t)-x'(t)\big) \big\| \, dt.
\end{eqnarray*}
According to (\ref{NVCondition2}) and (\ref{NVCondition3}),
the terms on $[0,T] \setminus K$ do not exceed
$$\frac{\delta}{3} \bigg( \|x(\cdot)-x'(\cdot)\|_\infty + \sum_{i=1}^d |\alpha_i - \alpha_i'| \bigg).$$
As in \cite{Ioffe}, pp. 246--251,
on the compact interval $K$ the third term can be estimated by
$$\frac{\delta}{3} \bigg( \|x(\cdot)-x'(\cdot)\|_\infty + \sum_{i=1}^d |\alpha_i - \alpha_i'| \bigg).$$
Finally, by (\ref{NVCondition1}), that last term does not exceed $\delta/3 \cdot \|x(\cdot)-x'(\cdot)\|_\infty$.
The proof is completed. \hfill $\square$ 

\begin{lemma} \label{LemmaNVUH2}
There exist $\sigma>0$ and $\Delta \in (0,1/d]$,
such that the inequality
\begin{equation} \label{NVUH2}
\Phi_2\big(x(\cdot),\alpha\big) - \Lambda_2\big(x(\cdot),\alpha\big) \leq \delta \sum_{i=1}^d \alpha_i
\end{equation}
holds for all $x(\cdot) \in V(\sigma)$ and all $\alpha \in \Sigma(\Delta)$.
\end{lemma}

\textbf{Proof.} The left-hand side of (\ref{NVUH2}) is equal to
$$\sum_{i=1}^{d} \int_0^\infty \Big(\chi_{M_i(\alpha_i)}(t) -\alpha_i \Big) \cdot
                \Big( \omega(t) \big[f \big(t,x(t),u_i(t)\big) - f\big(t,x(t),u_*(t)\big)\big]\Big) \, dt.$$
We take $M_i(\alpha_i) \subseteq K$ into account.
By the choice of $K$, it follows from (\ref{NVCondition4})
$$\sum_{i=1}^{d} \int_{[0,T] \setminus K} \alpha_i \cdot
                 \omega(t) \Big(f \big(t,x(t),u_i(t)\big) - f\big(t,x(t),u_*(t)\big)\Big) \, dt
    \leq \frac{\delta}{2} \sum_{i=1}^d \alpha_i.$$
Furthermore, the number $T>0$ is chosen,
such that $u_i(t)=u_*(t)$ for all $t>T$ and $i=1,...,d$.
This yields
$$\sum_{i=1}^{d} \int_T^\infty \alpha_i \cdot \Big( \omega(t) \big[f \big(t,x(t),u_i(t)\big) - f\big(t,x(t),u_*(t)\big)\big]\Big) \, dt = 0.$$
As in \cite{Ioffe}, pp. 246--251,
on the compact interval $K$ we obtain 
$$\sum_{i=1}^{d} \int_K \Big(\chi_{M_i(\alpha_i)}(t) -\alpha_i \Big) \cdot
                \Big( \omega(t) \big[f \big(t,x(t),u_i(t)\big) - f\big(t,x(t),u_*(t)\big)\big]\Big) \, dt
  \leq \frac{\delta}{2} \sum_{i=1}^d \alpha_i.$$
The proof is completed. \hfill $\square$
\subsection{Linear Differential Equations}
The following study of linear differential equations of the form
\begin{equation} \label{LinearDE}
\dot{x}(t)= A(t)x(t)+a(t)
\end{equation}
is closely related to the finite horizon case in \cite{Ioffe}.
Analogously to \cite{Ioffe}, $t \to A(t)$ is a mapping of $\R_+$ into the space of linear operators from $\R^n$ to $\R^n$,
and $a(t):\R_+ \to \R^n$ is a vector-valued function.

\begin{lemma} \label{LemmaLinearDE1}
Assume that the mapping $t \to A(t)$ and the vector-valued function $a(t)$ are summable on $\R_+$.
Then, for every vector-valued function $z(\cdot) \in C_{\lim}(\R_+,\R^n)$ and every $\tau \in [0,\infty]$,
there exists a unique vector-valued function $x(\cdot) \in C_{\lim}(\R_+,\R^n)$ such that
$$x(t)=z(t) + \int_{\tau}^t [A(s) x(s)+a(s)] \, ds$$
for all $t \in [0,\infty]$.
\end{lemma}

\textbf{Proof.} According to the assumptions on the mappings $t \to A(t), t \to a(t)$,
the operator $Q$,
$$x(\cdot) \to [Qx(\cdot)](t)=z(t) + \int_{\tau}^t [A(s) x(s)+a(s)] \, ds,$$
maps the space $C_{\lim}(\R_+,\R^n)$ into itself.
Then Lemma \ref{LemmaLinearDE1} is shown,
if the fixed-point equation $x(\cdot) = Qx(\cdot)$ possesses a unique solution.
We set
$$c(t) = \|A(t)\|, \qquad C(t) = \int_\tau^t c(s) \, ds, \qquad c_0 = \int_0^\infty c(s) \, ds.$$
As in \cite{Ioffe},
it follows the estimation
$$\big\|Q^m \big(x_1(\cdot) - x_2(\cdot)\big) \big\|_\infty \leq \frac{c_0^m}{m!} \cdot \| x_1(\cdot) - x_2(\cdot) \|_\infty.$$
The required result follows from the fixed-point theorem of Weissinger (cf. \cite{HeuserGD}). \hfill $\square$

\begin{corollary} \label{CorollaryDE}
The linear operator $F_x\big(x_*(\cdot),u_*(\cdot)\big)$ maps the space $C_{\lim}(\R_+,\R^n)$ onto itself.
\end{corollary}

\begin{lemma} \label{LemmaLinearDE2}
Let the mapping $t \to A(t)$ and the vector-valued function $a(t)$ satisfy all conditions in Lemma \ref{LemmaLinearDE1}.
Then, for every $z \in \R^n$ and every $\tau \in [0,\infty]$,
there exists a unique solution $x(t)$ of (\ref{LinearDE}) on $[0,\infty]$ such that $x(\tau)=z$.
\end{lemma}
\subsection{Convex Analysis}
State constraints on the infinite horizon is a special and rarely considered topic.
At this point, we would recall the discussion in Section \ref{SectionSpace}.
For the present particular mappings 
we state the relation between the one-sided directional derivative considered in Ioffe \& Tichomirov \cite{Ioffe} and
Clarke's generalized gradient stated in \cite{Clarke},
$$f'(x_0;x)= \lim_{\lambda \to 0^+} \frac{f(x_0+\lambda x)-f(x_0)}{\lambda}, \quad
  f^\circ(x;z) = \limsup_{\substack{y \to x \\ \lambda \to 0^+}} \frac{f(y + \lambda z) - f(y)}{\lambda}.$$

\begin{lemma}
Let $X$ and $Y$ be Banach spaces,
let $G:X \to Y$ be Fr\'echet differentiable at $x_0$,
and let the proper convex function $\gamma:Y \to \R$ be continuous at $G(x_0)$.
Then, the function $f(x)=\gamma\big(G(x)\big)$ is regular in $x_0$, e.g.
$$f'(x_0;x)=f^{\circ}(x_0;x) \qquad \mbox{for all } x \in X.$$
\end{lemma}

Now, with the mappings $G_j$ defined in Section \ref{SectionEP} we consider the following functions on $C_{\lim}(\R_+,\R^n)$:
$$f_j\big(x(\cdot)\big) = \gamma\big( G_j\big(x(\cdot)\big)\big) = \max_{t \in \overline{\R}_+} g_j\big(t,x(t)\big), \quad j=1,...,l.$$

\begin{lemma} \label{LemmaZB2}
The function $\gamma:C_{\lim}(\R_+,\R) \to \R$ is proper, convex and continuous.
Moreover, by Lemma \ref{LemmaZB1} the mappings $G_j:C_{\lim}(\R_+,\R^n) \to C_{\lim}(\R_+,\R)$ are Fr\'echet differentiable at the point
$x_*(\cdot) \in \mathscr{X}_{Lip}$.
Consequently, the functions $f_j(x)=\gamma\big(G_j(x)\big)$, $j=1,...,l,$ are regularly locally convex in $x_*(\cdot)$ (cf. \cite{Ioffe}).
\end{lemma}

Finally, the standard arguments deliver the following representation (cf. \cite{Ioffe,DuboMil,Clarke}):

\begin{lemma} \label{LemmaZBSubdifferential}
For any $j=1,...,l$ the subdifferential of the function $f_j$ at the point $x_*(\cdot)$
contains those and only those continuous linear functions $x^* \in C_{\lim}^*(\R_+,\R^n)$
that can be represented in the form
$$\langle x^* , x(\cdot) \rangle = \int_0^\infty \langle g_{jx}\big(t,x_*(t)\big),x(t) \rangle \, d\mu_j(t),$$
where $\mu_j \in \mathscr{M}(\overline{\R}_+)$ is non-negative, having norm $1$ and supported on the set
$$T_j=\{t \in \overline{\R}_+ \,|\, g_j\big(t,x_*(t)\big) = f_j\big(x(\cdot)\big) \}.$$
\end{lemma}
\subsection{Completion of the Proof}
The optimal control problem (\ref{Aufgabe1})--(\ref{Aufgabe5}) in form of the extremal problem
(\ref{Extremalaufgabe1})
satisfies the requirements of Theorem \ref{TheoremLMR}.
The assumptions are shown in the previous parts of the proof
(more precisely: (a$_1$) by Lemma \ref{LemmaDynamik2}, (a$_2$) by Lemma \ref{LemmaZielfunktional},
(a$_3$) by Lemma \ref{LemmaZB2}, (b) by Lemma \ref{LemmaNVUH1} and Lemma \ref{LemmaNVUH2},
and (c) by Lemma \ref{CorollaryDE}). \\
We obtain from (\ref{SatzLMR1}) that for all $x(\cdot) \in C_{\lim}(\R_+,\R^n)$, the following variational equality must hold: 
\begin{eqnarray}
0 &=& \lambda_0 \cdot \int_0^\infty \omega(t) \langle f_x\big(t,x_*(t),u_*(t)\big),x(t) \rangle \, dt \nonumber \\
  & & + \int_0^\infty \bigg[ x(t)-x(0) - \int_0^t \varphi_x\big(s,x_*(s),u_*(s)\big) x(s) \,ds \bigg]^T d\nu(t) \nonumber \\
  & & + \langle l_0 , h_0'\big(x_*(0)\big) x(0) \rangle
      + \lim_{t \to \infty}\langle l_1 , h_{1x}\big(t,x_*(t)\big) x(t) \rangle \nonumber \\
  & & \label{Beweisschluss1} + \sum_{j=1}^l \lambda_j \int_0^\infty g_{jx}\big(s,x_*(s)\big) \, d\tilde{\mu}_j(t).
\end{eqnarray}
According to Riesz' representation theorem
$\nu=(\nu_1,...,\nu_n)$ denotes a vector of Borel measures $\nu_1,...,\nu_n \in \mathscr{M}(\overline{\R}_+)$.
Moreover, one can assume without loss of generality (cf. \cite{Ioffe}),
all Borel measures $\mu_j=\lambda_j\tilde{\mu}_j \in \mathscr{M}(\overline{\R}_+)$ are non-negative and supported on the sets
$$T_j=\big\{ t \in \overline{\R}_+ \,|\, g_j\big(t,x_*(t)\big)=0\big\}.$$
In equation (\ref{Beweisschluss1}) any term is absolutely integrable.
Due to the Fubini theorem we change the order of integration in the second term and bring the relation written above to the form
\begin{eqnarray}
0 &=&  \int_0^\infty \bigg[ \lambda_0 \omega(t)f_x\big(t,x_*(t),u_*(t)\big)
             - \varphi^T_x\big(t,x_*(t),u_*(t)\big) \int_t^\infty d\nu(s) \bigg]^T x(t) dt \nonumber \\
  & & + \int_0^\infty \langle x(t),d\nu(t) \rangle
      + \Big\langle {h_0'}^T\big(x_*(0)\big)\, l_0 -\int_0^\infty d\nu(t), x(0) \Big\rangle
      + \lim_{t \to \infty} \langle h_{1x}^T\big(t,x_*(t)\big)\, l_1 , x(t) \rangle \nonumber \\
  & & \label{Beweisschluss2} + \sum_{j=1}^l \int_0^\infty g_{jx}\big(s,x_*(s)\big) \, d\mu_j(t).
\end{eqnarray}
The right-hand side of (\ref{Beweisschluss2}) defines a linear and continuous functional on the space $C_{\lim}(\R_+,\R^n)$.
Applying the Riesz representation theorem and denoting $p(t)=\displaystyle \int_t^\infty \, d\nu(s)$, we obtain
\begin{eqnarray*}
p(t) &=& - \lim_{t \to \infty}h_{1x}^T\big(t,x_*(t)\big) l_1 + \int_t^\infty H_x\big(s,x_*(s),u_*(s),p(s),\lambda_0\big) \, ds
         - \sum_{j=1}^l \int_t^\infty g_{jx}\big(s,x_*(s)\big) \, d\mu_j(s), \\
p(0) &=& {h_0'}^T\big(x_*(0)\big)l_0.
\end{eqnarray*}
Thus, (\ref{PMP1}) and (\ref{PMP2}) are proven.
The relation (\ref{SatzLMR2}) is equivalent to
$$\int_0^\infty H\big(t,x_*(t),u_*(t),p(t),\lambda_0\big) \, dt
  \geq \min_{u(\cdot) \in \mathscr{U}} \int_0^\infty H\big(t,x_*(t),u(t),p(t),\lambda_0\big) \, dt.$$
Finally, the pointwise maximum condition (\ref{PMP3}) follows by standard arguments for the Lebesgue points (cf. \cite{Natanson}).
The proof of the Pontryagin maximum principle is completed.
\section{Transversality Conditions and the Normal Form in the Absence of State Constraints} \label{SectionTransversality}
In the following we differ between special types of boundary conditions at infinity:

\begin{definition}
We say (\ref{Aufgabe1})--(\ref{Aufgabe4}) is a problem with
\begin{enumerate}
\item[(1)] free endpoint iff $\displaystyle\lim_{t \to \infty} x(t)$ is free;
\item[(2)] fixed  endpoint iff $\displaystyle\lim_{t \to \infty} x(t) =x_1$ with a fixed vector $x_1 \in \R^n$;
\item[(3)] mixed, partly free and partly fixed, endpoint iff $\displaystyle\lim_{t \to \infty} x^{i}(t)$ is free for $i=1,...,l$ and
           $\displaystyle\lim_{t \to \infty} x^{i}(t) =x^{i}_1$ with fixed numbers $x^{i}_1$ for $i=l+1,...,n$.
\end{enumerate}
For convenience reasons, $\displaystyle\lim_{t \to \infty} x(t) =x_1$ denotes the mixed endpoint.
\end{definition}

According to the transversality conditions in Theorem \ref{SatzPMP} we obtain immediately the following versions of
they so called natural transversality conditions:

\begin{lemma} \label{LemmaNatural}
Let $p(\cdot)$ be the vector-valued function stated in Theorem \ref{SatzPMP}.
Consider the problem (\ref{Aufgabe1})--(\ref{Aufgabe4}) with
\begin{enumerate}
\item[(a)] free endpoint, then the following conditions hold:
           $$\lim_{t \to \infty} p(t) =0, \qquad \lim_{t \to \infty} \langle p(t),x(t) \rangle = 0 \quad \mbox{ for all } x(\cdot) \in W^1_\infty(\R_+,\R^n);$$
\item[(b)] mixed endpoint, then the following condition holds:
           $$\lim_{t \to \infty} \langle p(t),x(t)-x_*(t) \rangle = 0 \quad \mbox{ for all admissible } x(\cdot) \in W^1_\infty(\R_+,\R^n).$$
\end{enumerate}
\end{lemma}

According to $\big(x_*(\cdot),u_*(\cdot)\big) \in \mathscr{A}_{\lim}$ the condition
\begin{equation} \label{BedingungNormalform}
\int_0^\infty \big\|\varphi_x\big(t,x_*(t),u_*(t)\big)\big\| \, dt < \infty
\end{equation}
is satisfied.
It follows by Lemma \ref{LemmaLinearDE2}, that the equation
$$\dot{z}(t) = -\varphi_x^T\big(t,x_*(t),u_*(t)\big) z(t), \qquad \lim_{t \to \infty} z(t)=0,$$
possesses the unique solution $z(t) \equiv 0$ on $\R_+$.
Hence we can state:

\begin{lemma} \label{LemmaNF}
In the problem (\ref{Aufgabe1})--(\ref{Aufgabe4}) with free endpoint the normal form of the maximum principle holds,
e.g. Theorem \ref{SatzPMP} is satisfied with $\lambda_0=1$.
\end{lemma}

Let $Y_*(t), Z_*(t)$ be the fundamental matrix solutions (normalized at $t=0$) of the homogeneous linear systems
$$\dot{y}(t)=\varphi_x\big(t,x_*(t),u_*(t)\big) y(t), \quad \dot{z}(t)=-\varphi_x^T\big(t,x_*(t),u_*(t)\big) z(t).$$
We consider the function $y(t)= Y_*(t) [Y_*(T)]^{-1} \xi$ with arbitrarily fixed $T \in \R_+$ and $\xi \in \R^n$.
According to (\ref{BedingungNormalform}), the function $y(\cdot)$ belongs to the space $W^1_\infty(\R_+,\R^n)$.
Furthermore,
in the problem (\ref{Aufgabe1})--(\ref{Aufgabe4}) with free endpoint the following formula holds (cf. \cite{AseKry}):
$$p(t)=Z_*(t)\bigg(Z^{-1}_*(T) p(T) + \int_T^t \big[Z_*(s)\big]^{-1} \omega(s) f_x\big(s,x_*(s),u_*(s)\big) \, ds \bigg), \quad T \in \R_+.$$
Consider $\langle p(t), y(t) \rangle$. Using $Z^{-1}_*(t)=Y^T_*(t)$ we find
$$\langle p(t)\,,\, y(t) \rangle = \Big\langle  p(T) - Z_*(T)\int_T^t \big[Z_*(s)\big]^{-1} \omega(s) f_x\big(s,x_*(s),u_*(s)\big) ds \,,\, \xi \big\rangle.$$
Passing to the limit $t \to \infty$ and taking into account Lemma \ref{LemmaNatural} and Lemma \ref{LemmaNF},
we obtain the representation formula of the adjoint $p(\cdot)$ in integral form:

\begin{lemma} \label{LemmaNormalform}
Let $p(\cdot)$ be the vector-valued function stated in Theorem \ref{SatzPMP}.
In the problem (\ref{Aufgabe1})--(\ref{Aufgabe4}) with free endpoint the adjoint $p(\cdot)$ possesses the representation
\begin{equation} \label{Integraldarstellung}
p(t)= Z_*(t) \int_t^\infty \big[Z_*(s)\big]^{-1} \omega(s) f_x\big(s,x_*(s),u_*(s)\big) \, ds,
\end{equation}
where $Z_*(t)$ is the fundamental matrix solution (normalized at $t=0$) of the linear system
$\dot{z}(t)=-\varphi^T_x\big(t,x_*(t),u_*(t)\big) z(t).$
\end{lemma}
Note that the representation formula (\ref{Integraldarstellung}) is the same as in \cite{AseKry,AseVel,AseVel2,AseVel3}.
In contrast to the theory in \cite{AseKry,AseVel,AseVel2,AseVel3} the Lemma \ref{LemmaNormalform} characterizes
strong local minimizer in the present infinite horizon optimal control problems. \\[2mm]
Now, let the density function $\omega(\cdot) \in L_1(\R_+,\R)$ possesses a limit at infinity,
e.g. $\displaystyle\lim_{t \to \infty} \omega(t) = 0$.
Then, in the problem (\ref{Aufgabe1})--(\ref{Aufgabe4}) with free endpoint,
the limits
$$\lim_{t \to \infty} \omega(t) f\big(t,x_*(t),u_*(t)\big)=0, \qquad
  \lim_{t \to \infty} \big\langle p(t)\,,\, \varphi\big(t,x_*(t),u_*(t)\big)\big\rangle=0$$
hold.
Consequently, the transversality condition of Michel (cf. \cite{Michel}) is satisfied:

\begin{lemma}
In the problem (\ref{Aufgabe1})--(\ref{Aufgabe4}) with free endpoint let $\displaystyle\lim_{t \to \infty} \omega(t) = 0$.
Then it follows
$$\lim_{t \to \infty} H\big(t,x_*(t),u_*(t),p(t),1\big)=0.$$
\end{lemma}

\section{Sufficiency Conditions} \label{SectionArrow}
In this Section,
we discuss Arrow type sufficiency conditions in the problem
\begin{eqnarray}
&& \label{AufgabeHB1} J\big(x(\cdot),u(\cdot)\big) = \int_0^\infty \omega(t) f\big(t,x(t),u(t)\big) \, dt \to \inf, \\
&& \label{AufgabeHB2} \dot{x}(t) = \varphi\big(t,x(t),u(t)\big), \\
&& \label{AufgabeHB3} x(0)=x_0, \qquad \lim_{t \to \infty} x(t)=x_1,\\
&& \label{AufgabeHB4} u(t) \in U \subseteq \R^m, \quad U \not= \emptyset, \\
&& \label{AufgabeHB5} g_j\big(t,x(t)\big) \leq 0, \quad t \in \R_+, \quad j=1,...,l.
\end{eqnarray}
To state the sufficiency conditions we introduce the Hamiltonian
$$\mathscr{H}(t,x,p) = \sup_{u \in U} H(t,x,u,p,1).$$

\begin{theorem} \label{SatzArrow} 
Let $\big(x_*(\cdot),u_*(\cdot)\big) \in \mathscr{A}_{adm} \cap \mathscr{A}_{Lip}$.
Moreover,
suppose the pair $\big(x_*(\cdot),u_*(\cdot)\big)$, together with an adjoint $p(\cdot)$,
satisfies the conditions (\ref{PMP1})--(\ref{PMP3}) of the Pontryagin maximum principle in normal form.
In addition, suppose that the Hamiltonian $\mathscr{H}(t,x,p)$ is a concave function of the variable $x$ on $V_\gamma$.
Then, the admissible pair $\big(x_*(\cdot),u_*(\cdot)\big)$ is a strong local minimizer
in problem (\ref{AufgabeHB1})--(\ref{AufgabeHB4}).
\end{theorem}

\begin{remark}
In contrast to the Pontryagin maximum principle,
the statement of the Arrow sufficiency conditions replaces $\big(x_*(\cdot),u_*(\cdot)\big) \in \mathscr{A}_{\lim}$ by
assumptions on concavity.
\end{remark}

\textbf{Proof.} 
In the proof we follow \cite{AseKry} and \cite{Seierstad}.
As shown in \cite{AseKry,Seierstad}, the difference
$$\Delta(T) = \int_0^T \omega(t)\big[f\big(t,x(t),u(t)\big)-f\big(t,x_*(t),u_*(t)\big)\big] \, dt$$
leads to the inequality
$$\Delta(T) \geq \langle p(T),x(T)-x_*(T)\rangle-\langle p(0),x(0)-x_*(0)\rangle$$
for any admissible pair $\big(x(\cdot),u(\cdot)\big)$ with $\|x(\cdot)-x_*(\cdot)\|_\infty \leq \gamma$.
Applying the natural transversality conditions in Lemma \ref{LemmaNatural},
we conclude
$$\lim_{T \to \infty} \Delta(T) \geq \lim_{T \to \infty} \langle p(T),x(T)-x_*(T)\rangle-\langle p(0),x(0)-x_*(0)\rangle=0$$
for any admissible pair $\big(x(\cdot),u(\cdot)\big)$ with $\|x(\cdot)-x_*(\cdot)\|_\infty \leq \gamma$.  \hfill $\square$

\begin{example}\label{ExampleRegulator1} We consider the linear-quadratic problem
\begin{eqnarray*}
&& \int_0^\infty \frac{1}{2}e^{-2t}\big(x^2(t)+u^2(t)\big) \, dt \to \inf, \\
&& \dot{x}(t)=2x(t)+u(t), \quad x(0)=2, \quad u(t) \in \R.
\end{eqnarray*}
Any admissible process $\big(x(\cdot),u(\cdot)\big)$ violates the requirement
$\big(x(\cdot),u(\cdot)\big) \in \mathscr{A}_{\lim}$ of Theorem \ref{SatzPMP}.
But we may apply the Arrow sufficiency conditions.
The conditions (\ref{PMP1})--(\ref{PMP3}) in normal form deliver the pair
$$x_*(t)=2e^{(1-\sqrt{2})t}, \quad u_*(t)= -2(1+\sqrt{2})e^{(1-\sqrt{2})t}$$
together with the adjoint $p(t)=u_*(t)e^{-2t}$.
Since the Hamiltonian is a concave function of the variable $x$,
the process $\big(x_*(\cdot),u_*(\cdot)\big)$ is a strong local minimizer.
\end{example}

\begin{definition}
We say the function $p(\cdot):\R_+ \to \R^n$ satisfies piecewise the adjoint equation (\ref{PMPZB1}) iff there
exist sequences $\{s_n\},\{\beta_{jn}\}$ and functions $\lambda_j(\cdot) \in L_1(\R_+,\R^n)$, $n=0,1,...,j=1,...,l$,
such that the following conditions hold:
\begin{enumerate}
\item $0=s_0<s_1<...$ and any finite interval $[0,T]$ containing only a finite number of elements $s_n$;
\item $\lambda_j(t)$ are non-negative and continuous in $(s_n,s_{n+1})$,  $n=0,1,...,\; j=1,...,l$;
\item $\displaystyle \dot{p}(t)= -H_x\big(t,x_*(t),u_*(t),p(t),1\big) + \sum_{j=1}^l \lambda_j(t) g_{jx}\big(t,x_*(t)\big),
       \quad t \in (s_n,s_{n+1})$;
\item $\displaystyle p(s_n^+) - p(s_n) = \sum_{j=1}^l \beta_{jn} g_{jx}\big(s_n,x_*(s_n)\big), \quad$
      $\displaystyle \sum_{n=0}^\infty \beta_{jn}<\infty,\quad$ $\beta_{jn} \geq 0,\quad$ $n=0,1,...,$ $j=1,...,l$.
      Herein, $p(s_n^+)$ denotes the right-sided limit.
\end{enumerate}
\end{definition}

\begin{theorem} \label{SatzArrowZB} 
Let $\big(x_*(\cdot),u_*(\cdot)\big) \in \mathscr{A}_{adm} \cap \mathscr{A}_{Lip}$.
Moreover,
suppose the pair $\big(x_*(\cdot),u_*(\cdot)\big)$, together with an adjoint $p(\cdot)$,
satisfies piecewise the conditions (\ref{PMPZB1})--(\ref{PMPZB3}) of Theorem \ref{SatzPMPZB} in normal form.
More precisely, there exist vectors $l_0,l_1 \in \R^n$ and numbers $\beta_j \geq 0$, $j=1,...,l$, such that
\begin{enumerate}
\item[(a)] the vector-valued function $p(\cdot):\R_+ \to \R^n$ satisfies piecewise the adjoint equation (\ref{PMPZB1}),
           \begin{eqnarray*}
           && \dot{p}(t)= -H_x\big(t,x_*(t),u_*(t),p(t),1\big) + \sum_{j=1}^l \lambda_j(t) g_{jx}\big(t,x_*(t)\big), \quad t \in (s_n,s_{n+1}), \\
           &&  p(s_n^+) - p(s_n) = \sum_{j=1}^l \beta_{jn} g_{jx}\big(s_n,x_*(s_n)\big), \quad n=0,1,...,
           \end{eqnarray*}
           the transversality conditions (\ref{PMPZB2}) and  (\ref{Transversality}),
           $$p(0)=l_0, \qquad \lim_{t \to \infty} p(t)= -l_1 - \lim_{t \to \infty} \sum_{j=1}^l \beta_j g_{jx}\big(t,x_*(t)\big),$$
           and for all admissible $x(\cdot)$ the natural transversality conditions
           $$\langle l_0,x(0)-x_*(0)\rangle=0, \qquad \langle l_1, \lim_{t \to \infty} x(t)-x_*(t)\rangle = 0;$$
\item[(b)] for almost all $t \in \R_+$ the maximum condition (\ref{PMPZB3}) holds:    
           $$H\big(t,x_*(t),u_*(t),p(t),1\big) = \max_{u \in U} H\big(t,x_*(t),u,p(t),1\big)\;$$
\item[(c)] for $j=1,...,l$, $n=0,1,...$ and $t \in (s_n,s_{n+1})$ the complementary slackness conditions hold:
           $$\beta_{jn} g_j\big(s_n,x_*(s_n)\big)=0, \qquad
             \lim\limits_{t \to \infty} \beta_j g_j\big(t,x_*(t)\big) =0, \qquad
             \lambda_j(t) g_j\big(t,x_*(t)\big) =0.$$
\end{enumerate} 
Suppose that the Hamiltonian $\mathscr{H}(t,x,p)$ is a concave function of the variable $x$ on $V_\gamma$.
Furthermore, suppose $g_j(t,x)$ are convex functions of the variable $x$ on $V_\gamma$.
Then, the admissible pair $\big(x_*(\cdot),u_*(\cdot)\big)$ is a strong local minimizer in problem
(\ref{AufgabeHB1})--(\ref{AufgabeHB5}).
\end{theorem}

\textbf{Proof.} By standard results on convex functions, we obtain:
$$g_j(t,x) \geq g_j\big(t,x_*(t)\big) + \big\langle g_{jx}\big(t,x_*(t)\big) , x-x_*(t) \big\rangle.$$
For any admissible $x(\cdot)$, this inequality together with the complementary slackness conditions deliver
$$0 \leq -\lim_{t \to \infty} \sum_{j=1}^l \beta_j g_j\big(t,x(t)\big)
    \leq -\lim_{t \to \infty} \sum_{j=1}^l \beta_j \big\langle g_{jx}\big(t,x_*(t)\big) , x(t)-x_*(t) \big\rangle.$$
Let $T \in \R_+ \setminus \{s_0,s_1,...\}$.
As shown in \cite{Seierstad}, the difference
$$\Delta(T) = \int_0^T \omega(t)\big[f\big(t,x(t),u(t)\big)-f\big(t,x_*(t),u_*(t)\big)\big] \, dt$$
leads to the inequality
$$\Delta(T) \geq \langle p(T),x(T)-x_*(T)\rangle -\langle p(0),x(0)-x_*(0)\rangle$$
for any admissible pair $\big(x(\cdot),u(\cdot)\big)$ with $\|x(\cdot)-x_*(\cdot)\|_\infty \leq \gamma$. \\
According to $\displaystyle \sum_{n=0}^\infty \beta_{jn}<\infty$ we have $\lim\limits_{n \to \infty} \|p(s_n^+) - p(s_n)\|=0$.
Applying the transversality conditions stated in Theorem \ref{SatzArrow} and taking into account $\lim\limits_{n \to \infty} \|p(s_n^+) - p(s_n)\|=0$,
we obtain
\begin{eqnarray*}
\lim_{T \to \infty} \Delta(T) &\geq& \lim_{T \to \infty} \langle p(T),x(T)-x_*(T)\rangle-\langle p(0),x(0)-x_*(0)\rangle \\
&=& \lim_{T \to \infty} \Big[ -\langle l_1,x(T)-x_*(T)\rangle + \langle l_0,x(0)-x_*(0)\rangle \Big] \\
& & + \lim_{T \to \infty} \bigg[ -\sum_{j=1}^l \beta_j \big\langle g_{jx}\big(T,x_*(T)\big) , x(T)-x_*(T) \big\rangle\bigg]
    \geq 0
\end{eqnarray*}
for any admissible pair $\big(x(\cdot),u(\cdot)\big)$ with $\|x(\cdot)-x_*(\cdot)\|_\infty \leq \gamma$. \hfill $\square$

\begin{example}
We consider the problem of resource extraction with waste in \cite{Seierstad}:
\begin{eqnarray*}
&& J\big(x(\cdot),y(\cdot),u(\cdot)\big) =\int_0^\infty e^{-rt}\big[f\big(u(t)\big)-ay(t)-qu(t)\big] \, dt \to \sup, \\
&& \dot{x}(t) = -u(t),\quad \dot{y}(t)=cf\big(u(t)\big), \quad x(0)=x_0>0,\quad y(0)=y_0\geq 0, \\
&& x(t) \geq 0, \qquad u(t) \geq 0, \qquad r,a,q,c >0, \qquad r - ac>0.
\end{eqnarray*}
Assume that $f$ is twice continuously differentiable, $f \geq 0$, $f' >0$, $f''<0$. 
For an economical interpretation cf. \cite{Seierstad}.
In contrast to \cite{Seierstad},
we replaced the condition $\liminf\limits_{t \to \infty} x(t) \geq 0$ by the state constraint $x(t) \geq 0$ on $\R_+$. \\
In this problem the optimality conditions have the form:
\begin{enumerate}
\item[(a)] The Pontryagin function is
           $$H(t,x,y,u,p_1,p_2,1) = p_1 (-u)+p_2 cf(u) + e^{-rt}[f(u)-ay-qu].$$
\item[(b)] The adjoints satisfying the equations
           $$p_1(t)=\int_t^\infty \, d\mu(s), \qquad \dot{p}_2(t)=a e^{-rt} \Rightarrow p_2(t)=-\frac{a}{r}e^{-rt} + K,$$
           where the measure $\mu$ is concentrated on the set $T=\{t \in \overline{\R}_+ \,|\, x_*(t)=0\}$.
           Therefore, the function $p_1(\cdot)$ is non-negative and monotonically decreasing.
           Moreover, the transversality conditions deliver $K=0$.
\item[(c)] The maximum condition can be reduced to
           \begin{eqnarray*}
           \lefteqn{\Big[ -p_1(t) u_*(t) +c p_2(t) f\big(u_*(t)\big) + e^{-rt}[f\big(u_*(t)\big)-qu_*(t)]\Big]} \\
           && =\max_{u \geq 0} \Big[ -p_1(t) u +c p_2(t) f(u) + e^{-rt}[f(u)-qu]\Big].
           \end{eqnarray*}
           Let $d=(r-ac)/r$. Then the representation of $p_2(t)$ leads to:
           $$\Big[ d f\big(u_*(t)\big)e^{-rt}-u_*(t)\big(p_1(t)+ qe^{-rt}\big)\Big]
             =\max_{u \geq 0} \Big[ d f(u)e^{-rt}-u\big(p_1(t)+ qe^{-rt}\big)\Big].$$
\end{enumerate}
We introduce the function $g(u)= d f(u)e^{-rt}-u\big(p_1(t)+ qe^{-rt}\big)$.
This function is twice continuously differentiable with
$$g'(u)= \big(d f'(u)-q\big)e^{-rt}-p_1(t), \quad g''(u)=df''(u)e^{-rt}, \quad d=\frac{r-ac}{r}>0.$$
We discuss the following three cases:
\begin{enumerate}
\item[(A)] $df'(0)\leq q$: This case leads to $g'(0) \leq 0$ and we obtain
           $$u_*(t) \equiv 0, \quad x_*(t) \equiv x_0, \quad
             y_*(t)=y_0 + cf(0)t, \quad p_1(t) \equiv 0, \quad p_2(t)=-\frac{a}{r}e^{-rt}.$$
\item[(B)] $df'(0)> q$ and $p_1(0)=0$:
           By $g'(u)=\big(d f'(u)-q\big)e^{-rt}=0$ we obtain the optimal control $u_*(t)=u_0>0$ on $\R_+$,
           in contradiction to $x_*(t) \geq 0$ for all $t \in \R_+$.
\item[(C)] $df'(0)> q$ und $p_1(0)>0$:
           Since $p_1(0)>0$, the resource will be completely extracted.
           That means, there exists $t'>0$ with $x_*(t)>0$ for $t \in [0,t')$ and $x_*(t)=0$ for $t\geq t'$.
           Consequently, $u_*(t)=0$ for $t\geq t'$. 
           For $t\geq t'$ the adjoint $p_1(\cdot)$ is monotonically decreasing and satisfies
           $$p_1(t) = \big(df'(0)-q\big)e^{-rt'} \mbox{ for } t \leq t', \qquad
             p_1(t) = \big(df'(0)-q\big)e^{-rt} \mbox{ for } t \geq t'.$$
           According to $p_1(0) >0$ it follows $t'<\infty$.
           We show the uniqueness of $t'>0$:
           For $\tau \geq 0$ let the family $u_\tau(\cdot)$ be defined by 
           $$f'\big(u_\tau(t)\big)=\frac{1}{d}(q+[df'(0)-q]e^{r(t-\tau)}), \quad t \in [0,\tau],$$
           and $u_\tau(t)=0$ for $t \geq \tau$.
           Since $f'\big(u_\tau(\tau)\big)=f'(0)$, the functions $u_\tau(\cdot)$ are continuous.
           Moreover, the relation $f'\big(u_\tau(t)\big) < f'\big(u_s(t)\big)$, e.g. $u_\tau(t) > u_s(t)$, holds for all $t \in [0,\tau]$ and all $\tau>s$.
           Therefore, the function
           $$U(\tau):= \int_0^\infty u_\tau(t) \, dt$$
           is strictly increasing with $U(0)=0$.
           Then the number $t'$ satisfies the condition $U(t')=x_0$.
\end{enumerate}
In the cases (A) and (C) the adjoints $p_1(\cdot),p_2(\cdot)$ satisfying piecewise the adjoint equation.
Moreover, the problem is linear in the state variables $x,y$.
Thus, the Hamiltonian $\mathscr{H}$ is a concave function of $x,y$ and $g(t,x,y)=-x$ is convex in $x,y$.
Therefore, in the cases (A) and (C) the sufficiency conditions of Theorem \ref{SatzArrowZB} are satisfied.
\end{example}
\section{Optimal Control on Finite and Infinite Horizon} \label{SectionRelation}
In infinite horizon optimal control there exist several definitions of optimality and a corresponding number of methods to proof necessary optimality conditions. 
One of the most popular methods is the finite horizon approximation.
But the pathologies arising by this method are documented in literature (cf. \cite{AseKry,Halkin}).
As an example we discuss the linearized Ramsey model:
\begin{example}
We consider the problem
\begin{eqnarray*}
&& J\big(x(\cdot),u(\cdot)\big)=\int_0^\infty e^{-\varrho t}\big(1-u(t)\big)x(t) \, dt \to \sup, \\
&& \dot{x}(t)=u(t)x(t), \qquad x(0)=x_0, \qquad u \in [0,1], \qquad \varrho \in (0,1).
\end{eqnarray*}
For any $T>\tau = T + \ln (1-\varrho)/\varrho$ the process 
$$x^T_*(t)= \left\{ \begin{array}{ll} e^t,& t \in [0,\tau), \\[1mm] e^\tau,& t \in [\tau,T], \end{array} \right. \quad
  u^T_*(t)= \left\{ \begin{array}{ll} 1,& t \in [0,\tau), \\[1mm] 0,& t \in [\tau,T], \end{array} \right.$$
delivers the global maximizer in the problem
\begin{eqnarray*}
&& J_T\big(x(\cdot),u(\cdot)\big)=\int_0^T e^{-\varrho t}\big(1-u(t)\big)x(t) \, dt \to \sup, \\
 && \dot{x}(t)=u(t)x(t), \qquad x(0)=x_0, \qquad u \in [0,1], \qquad \varrho \in (0,1).
\end{eqnarray*}
But the limit $T \to \infty$ leads to the process
$$x_*(t) = e^t, \qquad u_*(t)=1, \qquad t \in \R_+,$$
which is the global minimizer in the problem on the infinite horizon.
\end{example}
Motivated by the pathologies, the question of the relation between finite horizon and infinite horizon optimal control 
problems arises. \\

One part of the answer can be given by the time transformation method:
Let $t:[0,1]\to \overline{\R}_+$ be a bijective function with $t(0)=0$.
Therefore, the mapping $t(\cdot)$ must satisfy the condition $t(1)=\infty$.
Then the time transformation (cf. \cite{Ioffe}) leads to the states $t(s)$, $y(s) = x\big(t(s)\big)$ and the control
$w(s) = u\big(t(s)\big)$,
and the problem (\ref{Aufgabe1})--(\ref{Aufgabe5}) possesses on $[0,1]$ the form
\begin{eqnarray*}
&& J\big(y(\cdot),w(\cdot)\big) = \int_0^1 v(s) \cdot \omega(s) f\big(t(s),y(s),w(s)\big) \, ds \to \inf, \\
&& y'(s) = v(s) \cdot \varphi\big(t(s),y(s),w(s)\big), \qquad t'(s) = v(s), \\
&& h_0\big(y(0)\big)=0,\quad h_1\big(t(1),y(1)\big)=0, \qquad t(0)=0, \quad t(1)=\infty, \\
&& w(s) \in U, \qquad v(s) > 0, \\
&& g_j\big(t(s),y(s)\big) \leq 0, \quad s \in [0,1], \quad j=1,...,l.
\end{eqnarray*}
This problem including the singularity $t(1)=\infty$,
which follows from the transformation of $\overline{\R}_+$ onto the finite interval $[0,1]$.
Or in other words, the infinite horizon possesses the nature of a singularity.
Therefore, the infinite horizon optimal control problem (\ref{Aufgabe1})--(\ref{Aufgabe5}) cannot be reduced to
a classical one. \\

On the other hand, consider the standard optimal control problem on the finite horizon $[t_0,t_1] \subset \R_+$:
\begin{eqnarray}
&& \label{AufgabeEH1} J\big(x(\cdot),u(\cdot)\big) = \int_{t_0}^{t_1} f\big(t,x(t),u(t)\big) \, dt \to \inf, \\
&& \label{AufgabeEH2} \dot{x}(t) = \varphi\big(t,x(t),u(t)\big), \\
&& \label{AufgabeEH3} h_0\big(x(0)\big)=0, \qquad h_1\big(x(t_1)\big)=0,\\
&& \label{AufgabeEH4} u(t) \in U \subseteq \R^m, \quad U \not= \emptyset, \\
&& \label{AufgabeEH5} g_j\big(t,x(t)\big) \leq 0, \quad t \in [t_0,t_1], \quad j=1,...,l.
\end{eqnarray}
In the problem (\ref{AufgabeEH1})--(\ref{AufgabeEH5}) let all the conditions of Section \ref{SectionProblem} be satisfied on
$$V_\gamma= \{ (t,x) \in [t_0,t_1] \times \R^n\,|\, \|x-x(t)\|\leq \gamma\}.$$
We introduce the mappings
$$\tilde{f}(t,x,u)= \left\{ \begin{array}{ll} f(t,x,u), & t \in [t_0,t_1], \\ 0, & t \not \in [t_0,t_1], \end{array} \right. \quad
  \tilde{\varphi}(t,x,u)= \left\{ \begin{array}{ll} \varphi(t,x,u), & t \in [t_0,t_1], \\ 0, & t \not \in [t_0,t_1], \end{array} \right. \quad
  \begin{array}{l} \tilde{h}_0(x)= h_0(x), \\ \tilde{h}_1(t,x)=h_1(x). \end{array}$$
Then the restrictions (\ref{condition1}), (\ref{condition2}) are satisfied.
Moreover, we define the functions
$$\tilde{g}_j(t,x) = \left\{ \begin{array}{ll} g_j(t_0,x) - (1-e^{(t-t_0)^2}), & t < t_0, \\
                                               g_j(t,x), & t \in [t_0,t_1], \\
                                               g_j(t_1,x) - (1-e^{(t-t_1)^2}), & t > t_1. \end{array}\right.$$
With the density function $\omega(t) = \chi_{[t_0,t_1](t)}$ the problem (\ref{AufgabeEH1})--(\ref{AufgabeEH5}) becomes the form
of an infinite horizon optimal control problem:
\begin{eqnarray}
&& \label{AufgabeUH1} J\big(x(\cdot),u(\cdot)\big) = \int_0^\infty \omega(t)\tilde{f}\big(t,x(t),u(t)\big) \, dt \to \inf, \\
&& \label{AufgabeUH2} \dot{x}(t) = \tilde{\varphi}\big(t,x(t),u(t)\big), \\
&& \label{AufgabeUH3} \tilde{h}_0\big(x(0)\big)=0, \qquad \lim_{t \to \infty} \tilde{h}_1\big(t,x(t)\big)=0,\\
&& \label{AufgabeUH4} u(t) \in U \subseteq \R^m, \quad U \not= \emptyset, \\
&& \label{AufgabeUH5} \tilde{g}_j\big(t,x(t)\big) \leq 0, \quad t \in \R_+, \quad j=1,...,l.
\end{eqnarray}

\begin{remark}
In $t=t_0$ and $t=t_1$ the mappings $\tilde{f}(t,x,u)$, $\tilde{\varphi}(t,x,u)$ may are discontinuous.
But the method of proof presented in this paper is still applicable.
\end{remark}

Then the pair $\big(\tilde{x}_*(\cdot),\tilde{u}_*(\cdot)\big)$ is a strong local minimizer in the problem
(\ref{AufgabeUH1})--(\ref{AufgabeUH5}) iff
the pair $\big(x_*(\cdot),u_*(\cdot)\big)$ with $\big(x_*(t),u_*(t)\big)= \big(\tilde{x}_*(t),\tilde{u}_*(t)\big)$ on 
$[t_0,t_1]$ is a strong local minimizer
in the problem (\ref{AufgabeEH1})--(\ref{AufgabeEH5}).
In particular, the sets
$$T_j = \big\{t \in [t_0,t_1] \,\big|\, g_j\big(t,x_*(t)\big)=0\big\}, \qquad 
  \tilde{T}_j = \big\{t \in \overline{\R}_+=[0,\infty] \,\big|\, \tilde{g}_j\big(t,\tilde{x}_*(t)\big)=0\big\}$$
coincide for any $j=1,...,l$.
Furthermore, the requirement $\big(x_*(\cdot),u_*(\cdot)\big) \in \mathscr{A}_{\lim}$ is satisfied
in the problem (\ref{AufgabeUH1})--(\ref{AufgabeUH5}).
Applying the Pontryagin maximum principle \ref{SatzPMP} in the problem (\ref{AufgabeUH1})--(\ref{AufgabeUH4})
and taking into account that $\dot{p}(t)=0$ for $t \not \in [t_0,t_1]$,
then there exist non-trivial multipliers $\lambda_0 \geq 0$, vectors $l_0 \in \R^{s_0}$, $l_1 \in \R^{s_1}$ and a vector-valued function
$p(\cdot):[t_0,t_1] \to \R^n$, such that
\begin{enumerate}
\item[(a)] the vector-valued function $p(\cdot)$ satisfies almost everywhere on $[t_0,t_1]$ the adjoint equation
           $$\dot{p}(t) = -\varphi_x^T\big(t,x_*(t),u_*(t)\big) p(t) + \lambda_0 f_x\big(t,x_*(t),u_*(t)\big)$$
           and the transversality conditions
           $$p(t_0) = {h_0'}^T\big(x_*(t_0)\big)l_0, \qquad p(t_1)= - {h_1'}^T\big(x_*(t_1)\big) l_1;$$
\item[(b)] for almost all $t \in [t_0,t_1]$ the maximum condition holds:    
           $$H\big(t,x_*(t),u_*(t),p(t),\lambda_0\big) = \max_{u \in U} H\big(t,x_*(t),u,p(t),\lambda_0\big).$$
\end{enumerate}
This is the maximum principle for the standard optimal control (\ref{AufgabeEH1})--(\ref{AufgabeEH4}) of Pontryagin et. al. in \cite{Pontrjagin}. \\
In the presence of state constraints the Pontryagin maximum principle \ref{SatzPMPZB} delivers the existence of
$\lambda_0 \geq 0$, $l_0 \in \R^{s_0}$, $l_1 \in \R^{s_1}$,
a vector-valued function $p(\cdot):[t_0,t_1] \to \R^n$, and non-negative regular Borel measures $\mu_j$, $j=1,...,l$,
on $[t_0,t_1]$ supported on the sets $T_j = \big\{t \in [t_0,t_1] \,\big|\, g_j\big(t,x_*(t)\big)=0\big\}$, respectively,
not all zero and such that
\begin{enumerate}
\item[(a)] the vector-valued function $p(\cdot)$ is a solution of the adjoint equation
           \begin{eqnarray*}
           p(t) &=& - {h_1'}^T\big(x_*(t_1)\big) l_1 + \int_t^{t_1} H_x\big(s,x_*(s),u_*(s),p(s),\lambda_0\big) \, ds,
                    - \sum_{j=1}^l \int_t^{t_1} g_{jx}\big(s,x_*(s)\big) \, d\mu_j(s) \\
           p(t_0) &=& {h_0'}^T\big(x_*(t_0)\big)l_0;                    
           \end{eqnarray*}
\item[(b)] for almost all $t \in [t_0,t_1]$ the maximum condition holds:    
           $$H\big(t,x_*(t),u_*(t),p(t),\lambda_0\big) = \max_{u \in U} H\big(t,x_*(t),u,p(t),\lambda_0\big).$$
\end{enumerate}
In addition, in $t=t_1$ the following transversality conditions hold:
\begin{eqnarray*}
p(t_1)= p(t_1-) &=& - {h_1'}^T\big(x_*(t_1)\big) l_1 - \sum_{j=1}^l g_{jx}\big(t_1,x_*(t_1)\big) \mu_j(\{t_1\}), \\
p(t_1+)-p(t_1) &=& \sum_{j=1}^l g_{jx}\big(t_1,x_*(t_1)\big) \mu_j(\{t_1\}).
\end{eqnarray*}
This is the complete statement of the Pontryagin maximum principle for standard optimal control problems with state constraints
(cf. \cite{Ioffe,DuboMil}).

\section{Conclusions}
In this paper we developed a novel approach to a class of infinite horizon optimal control problems.
The result is a complete set of necessary optimality conditions in form of the Pontryagin maximum principle.
Moreover, we have shown the validity of several transversality conditions and Arrow type sufficiency conditions.
In order to proof the maximum principle we introduced new elements
in the needle variation method on an unbounded time interval and with an arbitrary summable density function,
in the field of convex analysis in the framework of continuous functions converging at infinity
and for linear differential equations on the infinite horizon.
Finally, we have demonstrated
that the standard control theory on the finite horizon is a particular case of
infinite horizon optimal control problems with bounded processes. \\
This paper provides different perspectives:
The trajectories converging at infinity is a suitable framework in problems with steady states and in problems
with the turnpike phenomenon.
In particular, the space of continuous functions converging at infinity may be useful in non-smooth control problems on
unbounded intervals.
Furthermore,
the introduction of spaces of functions converging at infinity in Section \ref{SectionSpace} can be directly extended to the
cases of
sequence spaces $\ell_p$, Lebesgue spaces $L_p(\R_+,\R^n)$ and Sobolev spaces $W^1_p(\R_+,\R^n)$ with $1 \leq p < \infty$.
The control problems with bounded processes
playing an essential role in recent articles on infinite horizon optimal control problems
in the discrete time framework (cf. \cite{Blot,BachirBlot}) and in the dynamic programming (cf. \cite{Frankowska,Frankowska2}).
The demonstrated degenerations in the convex analysis in the space of continuous functions vanishing at infinity
carry over to the discrete time framework in the space of zero sequences.
An alternative will be the space of convergent sequences.
Moreover,
the achieved results are a solid base to discuss
the links between Pontryagin's maximum principle and Bellman's principle in infinite horizon optimal control
with bounded processes.

\end{document}